\documentclass[12]{article}
\usepackage{graphicx}
\usepackage[top=1in, bottom=1in, left=1in, right=1in]{geometry}
\newtheorem{theorem}{Theorem}
\usepackage{color}
\usepackage{amsmath}

\begin{document}
\title{On the Existence and Frequency Distribution of the Shell Primes}
\author{Michael P. May}
\date{August 13, 2015}
\maketitle
\begin{center}
\textit{In memory of Bill Joe Rodman Jr., brother, friend and mentor...}
\end{center}

\begin{abstract}
This research presents the results of a study on the existence and frequency distribution of the shell primes defined herein as prime numbers that result from the calculation of the "half-shell" of an p-dimensional entity of the form $n^p-(n-1)^p$ where power $p$ is prime and base $n$ is the realm of the positive integers.  Following the introduction of the shell primes, we will look at the results of a non-sieving application of the Euler zeta function to the prime shell function as well as to any integer-valued polynomial function in general which has the ability to produce prime numbers when power $p$ is prime.  One familiar with the Euler zeta function, which established the remarkable relationship between the prime and composite numbers, might naturally ponder the results of the application of this special function in cases where there is no known way to sieve composite numbers out of the product term in this famous equation. Such would be case when an infinite series of numbers to be analyzed are calculated by a polynomial expression that yields successively increasing positive integer values and which has as its input domain the positive integers themselves. In such cases there may not be an intuitive way to eliminate the composite terms from the product term in the Euler zeta function equation by either scaling a previous prime number calculation or by employing predictable values of the domain of the function which would render outputs of the polynomial prime. So the best one may be able to hope for in these cases is to calculate some value to be added or subtracted from unity in the numerator above the product term in the Euler Zeta function to make both sides of that equation equal with the expectation that that value could be used to predict the number of prime numbers that exist as outputs of the polynomial function for some limit less than or equal to x of the input domain.\\\\\\
\end{abstract}

\noindent This research introduces the results of a study on the existence and frequency distribution of the shell prime numbers.  Following the introduction of the shell primes, we will take a look at non-sieving applications of the Euler zeta function for integer-valued polynomials in general.\\\\We begin with a definition of the shell primes as they occur in nature.  The shell primes are defined herein as prime numbers that result from the calculation of the "half-shell" of an a-dimensional entity of the form

$$n^a-(n-1)^a$$

\noindent for positive integers $n\ge2$ and $a\ge2$.  This shell can intuitively be thought of as the partial outer boundary, or half-shell, of an a-dimensional entity of base $n$.  The first two occurrences are illustrated in Figs. 1 and 2.\\\\

\includegraphics[scale=0.029]{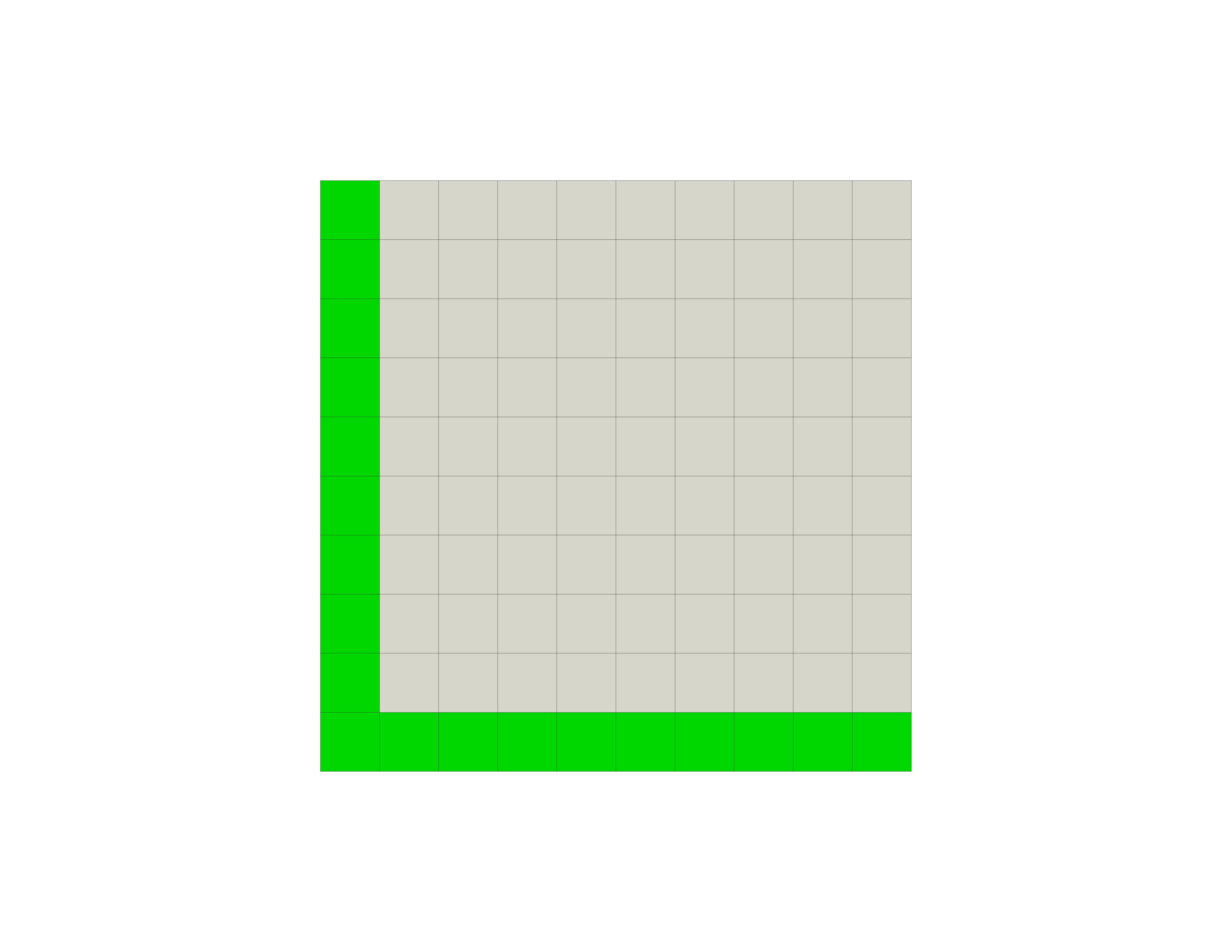}\includegraphics[scale=0.029]{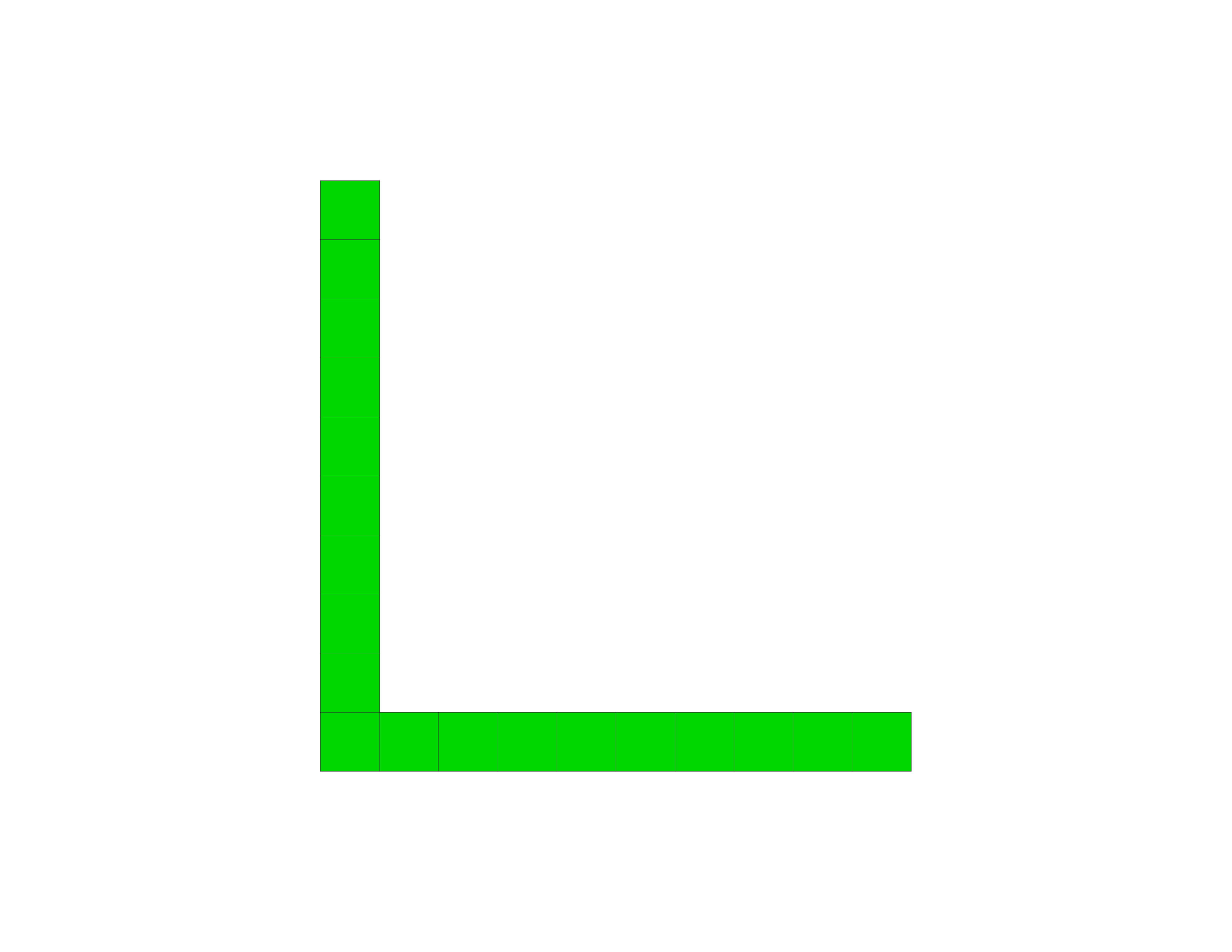}
\begin{center}
\textit{Fig. 1 - Half-shell of a 2-dimensional entity} 
\end{center}
\includegraphics[scale=0.03]{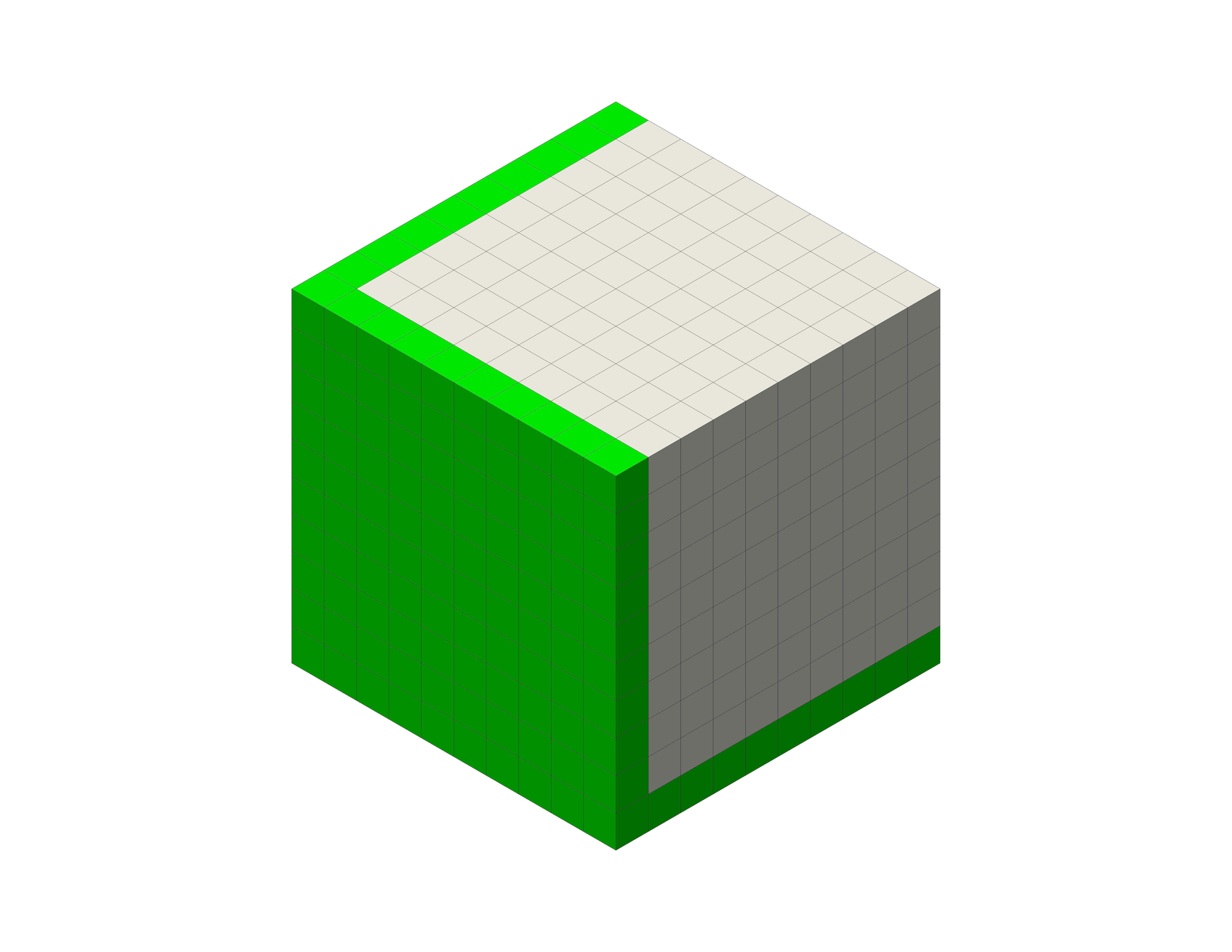}\includegraphics[scale=0.03]{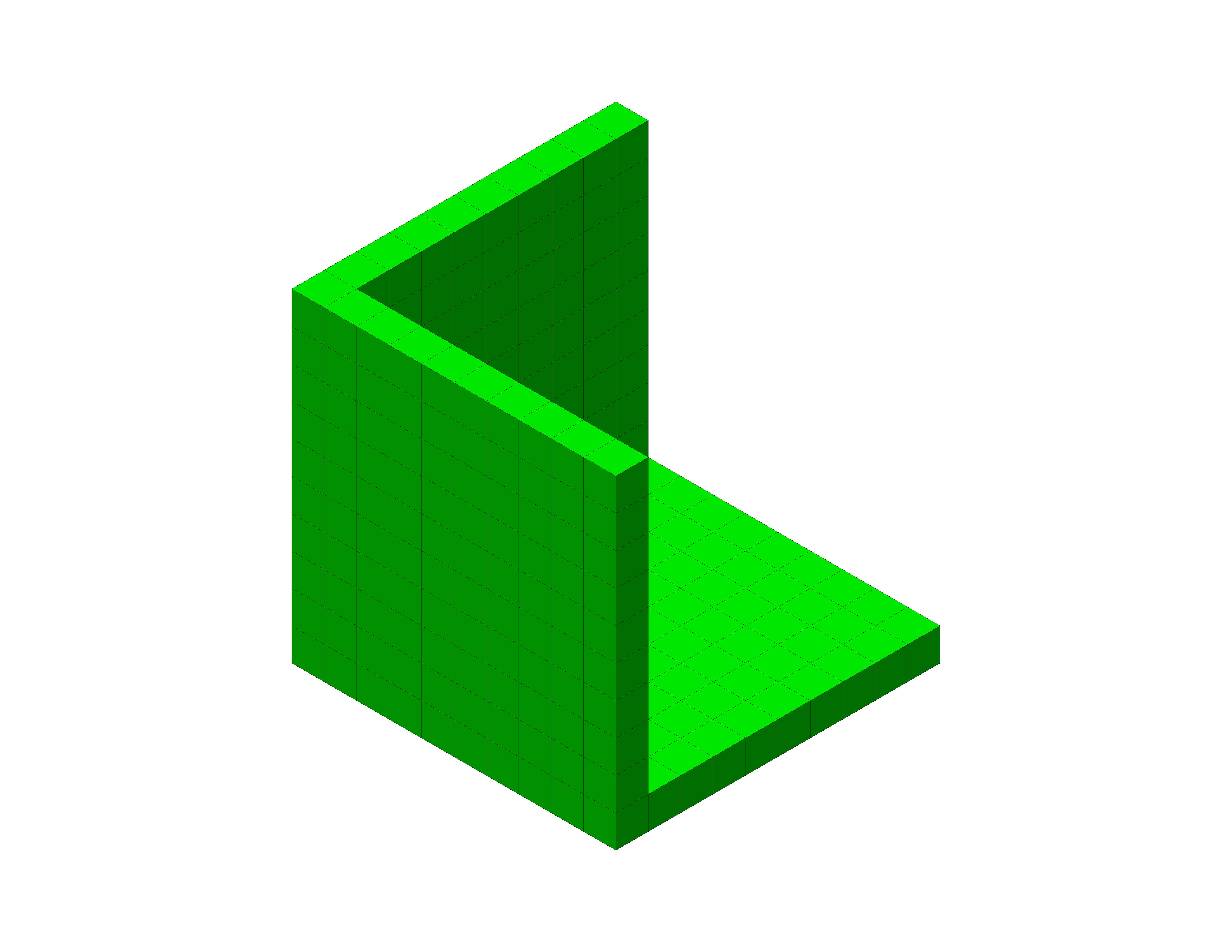}
\begin{center}
\textit{Fig. 2 - Half-shell of a 3-dimensional entity}
\end{center}

\noindent The shell of an $p$-dimensional entity is further defined by the function $n^p-(n-1)^p$ for any positive integer base $n\ge2$ and any positive prime power $p\ge2$.  We will refer to the function $n^p-(n-1)^p$ as the "prime shell" function since the base $n$ is taken to some prime power $p$.  The "shell primes" thus are prime numbers that result from the calculation of the prime shell function of base $n \ge 2$ for positive prime powers $p \ge 2$.\\

\noindent Within the framework of the above definitions, we now postulate the following theorem:\\

\begin{theorem}
Only prime shells, i.e., shells of the form $n^p-(n-1)^p$ for integer base $n \ge 2$ and prime power $p \ge 2$ will produce prime numbers.  Prime numbers do not exist for shells of the form $n^c-(n-1)^c$, where $c$ is a composite number, and this holds true regardless of whether the base $n$ is prime or composite.
\end{theorem}

\noindent In other words, shell primes are generated from prime shells only, or shells which are calculated using prime powers.  However, shell primes occur for both prime and composite values of the base $n$ when the power $p$ is prime, and they occur in a seemingly random fashion.  Prime numbers, however, do not exist for shells which are calculated using composite powers in the shell function $n^a-(n-1)^a$ regardless of whether the base $n$ is prime or composite.\\

\noindent The first few rows of the shell function are expanded in Fig. 3 using coefficients that are found in Paschal's triangle.  The prime shells, i.e., shells which will yield prime numbers, are illustrated in red.  Note that the leading coefficients of the prime shells are prime numbers themselves.  Also, note how all the coefficients of the powers of base $n$ in the shell functions in the prime rows are divisible by the first prime number in their row according to Fermat's Little Theorem.

\begin{center}
\color{red}
\item $2n-1$
\item $3n^2-3n+1$
\color{black}
\item $4n^3-6n^2+4n-1$
\color{red}
\item $5n^4-10n^3+10n^2-5n+1$
\color{black}
\item $6n^5-15n^4+20n^3-15n^2+6n-1$
\color{red}
\item $7n^6-21n^5+35n^4-35n^3+21n^2-7n+1$
\color{black}
\item $8n^7-28n^6-56n^5-70n^4+56n^3-28n^2+8n-1$
\color{black}
\item $9n^8-36n^7+84n^6-126n^5+126n^4-84n^3+36n^2-9n+1$
\color{black}
\item $10n^9-45n^9+120n^7-210n^6+252n^5-210n^4+120n^3-45n^2+1$
\color{red}
\item $11n^{10}-55n^9+165n^8-330n^7+462n^6-462n^5+330n^4-165n^3+55n^2-11n+1$
\color{black}
\item $12n^{11}-66n^{10}+220n^9-495n^8+792n^7-924n^6+792n^5-495n^4+220n^3-66n^2+2n-1$
\color{red}
\item $13n^{12}-78n^{11}+286n^{10}-715n^9+1287n^8-1716n^7+1716n^6-1287n^5+715n^4-286n^3+78n^2-13n+1$
\end{center}
\begin{center}
\bf.
\end{center}
\begin{center}
\bf.
\end{center}
\begin{center}
\bf.
\end{center}

\begin{center}
\textit{Fig. 3 - Expansion of the shell function $n^a-(n-1)^a$ for power $a\ge2$} 
\end{center}

\noindent It was observed that prime numbers generated by the prime shell function seem to arise in no less of a random fashion than do prime numbers on the real number line as can be seen in Figs. 4 and 5.\\

\includegraphics[scale=0.4]{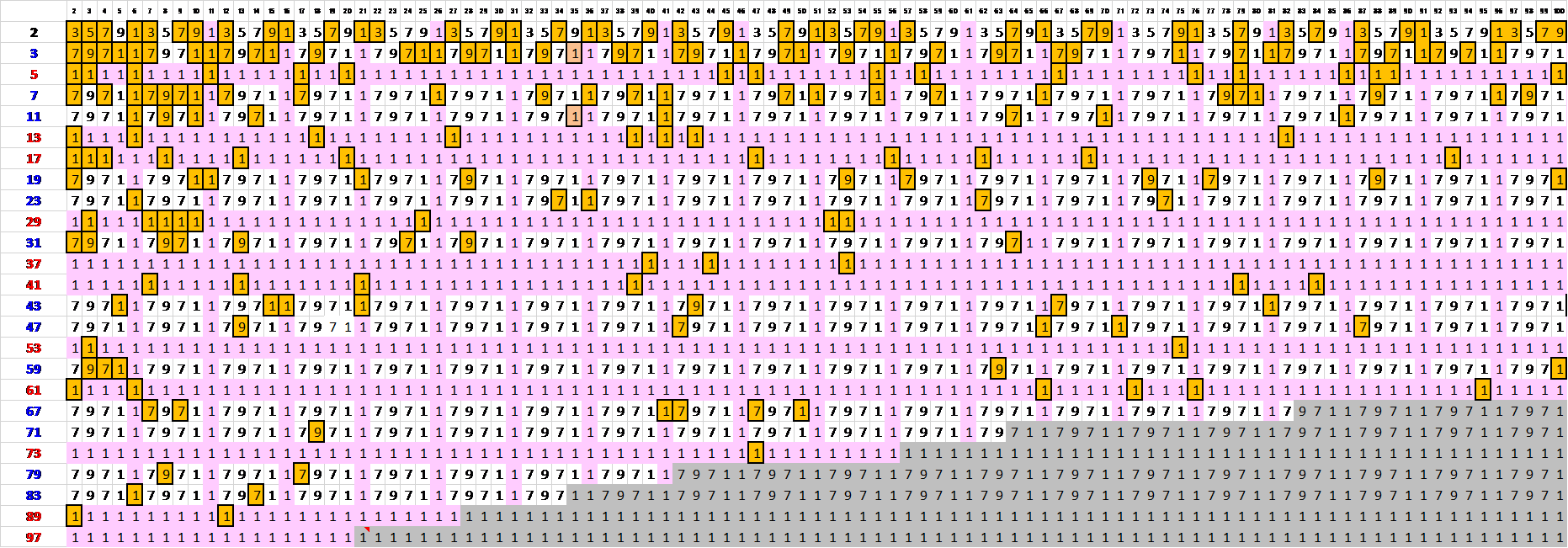}
\begin{center}
\textit{Fig. 4 - Spreadsheet matrix illustrating the distribution of the shell primes along the real number line calculated from the prime shells of the form $n^p-(n-1)^p$ for $2\le n \le100$ (columns) and for prime power $2\le p\le100$ (rows).  Descending rows correspond to increasing shell power $p$, and ascending columns to the right correspond to increasing base $n$.  Shells calculated with composite powers are omitted from the rows since they contain no prime numbers.  The grayed out sections denote the range of numbers having a number of digits beyond the capability of the author's resources to check for primality (128-digit limit).}\cite{TheNumberEmpire} 
\end{center}

\includegraphics[scale=0.4]{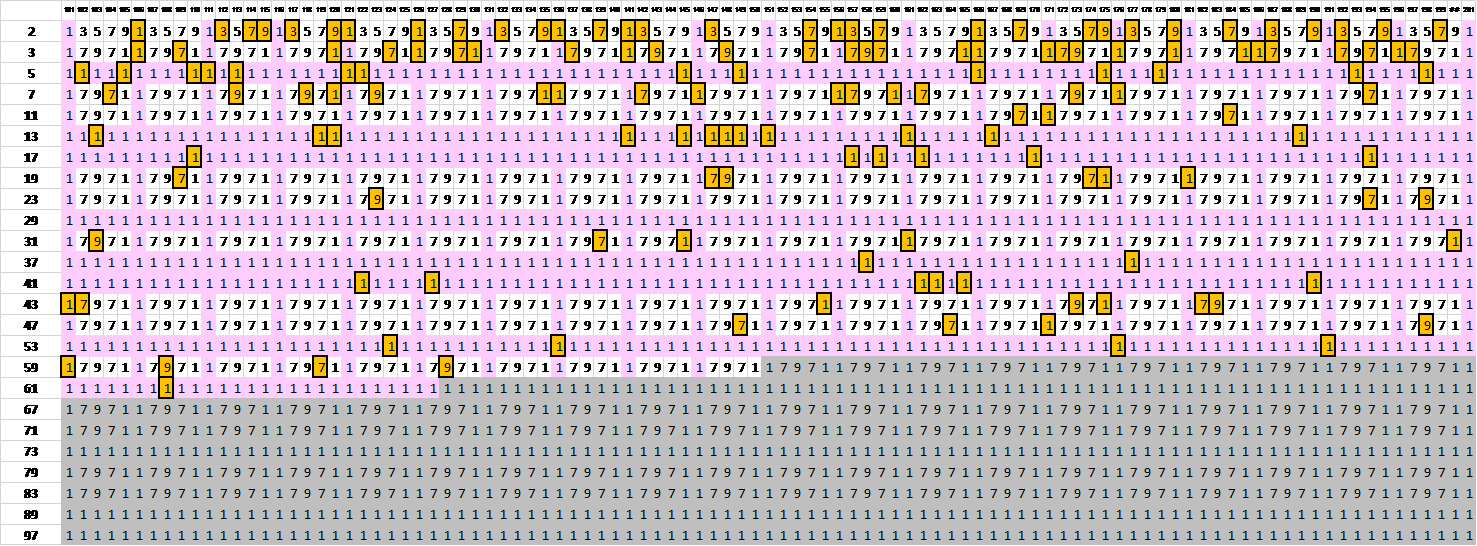}
\begin{center}
\textit{Fig. 5 - Continuation of the spreadsheet matrix depicted in Fig. 4 illustrating the distribution of the shell primes along the real number line calculated from the prime shells of the form $n^p-(n-1)^p$ for $100 < n \le 200$ (columns) and for prime power $2 \le p\le100$ (rows).  Descending rows correspond to increasing shell power $p$, and ascending columns to the right correspond to increasing base $n$.  Shells calculated with composite powers are omitted from the rows since they contain no prime numbers.  The grayed out sections denote the range of numbers having a number of digits beyond the capability of the author's resources to check for primality (128-digit limit).}\cite{TheNumberEmpire} 
\end{center}

\noindent One observation made from Figs. 4 and 5 is that as one looks downward and to the right in the charts, the prime numbers appear to be less frequent and further apart.  This is also evident by a few shell prime sequences that the author submitted to the Online Encyclopedia of Integer Sequences \cite{OnlineEncyclopediaofIntegerSequences} within the range of base $100 < n\le200$:\\

\color{blue} \underline {A254298} \color{black} \qquad Primes of the form $120^p-119^p$ for $p=2,3,7,13$:

239, 43841, 20386538221561, 110287289683553081554913641\\

\color{blue} \underline {A255387} \color{black} \qquad Primes of the form $141^p-140^p$ for $p=2,3,13$:

281, 59221, 769449701919846533025514621\\

\color{blue} \underline {A255388} \color{black} \qquad Primes of the form $157^p-156^p$ for $p=2,3,7,17$:

313, 73477, 102850464108757, 2202194587566133922938215539676032221\\

\color{blue} \underline {A255389} \color{black} \qquad Primes of the form $166^p-165^p$ for $p=2,3,5$:

331, 82171, 3751197451\\

\color{blue} \underline {A255390} \color{black} \qquad Primes of the form $173^p-172^p$ for $p=3,7,43$:

89269, 184438202074309, 37988946604016968734924614832155863...
\\

\noindent A second observation made from Figs. 4 and 5 is that all prime numbers generated from the prime shell function $n^p-(n-1)^p$ for prime power $p > 2$ end in either 1, 7, or 9.  Compare this result to the prime numbers that exist in the realm of the positive integers where all prime numbers beyond $2$ and $5$ end in either 1, 3, 7 or 9 and which, for some upper limit $ \le x$, are divided equally among those four ending digits \cite{TerenceTaoPrimeNumbers}.  This leads to a third theorem introduced as follows:

\begin{theorem}
All prime numbers calculated by the prime shell function $n^p-(n-1)^p$ for prime power $p > 2$ end in either 1, 7 or 9 with no primes ending in 3 resulting from those calculations.
\end{theorem}

\noindent Theorem 2 begs the question: Can prime numbers be calculated by other polynomial functions which would reduce the number of ending digits down to two (e.g., 1 and 7, or 7 and 9, or 1 and 9) or even down to one digit?\\

\noindent A third observation is in regard to the repeating pattern of the last digits in the calculations yielded by the general shell function $n^a-(n-1)^a$ for base $n \ge 2$ and power $a \ge 2$, prime and composite.  Please see Fig. 6 for a snapshot of the repeating pattern of the last digits in the shell function calculations for base $n \ge 2$ and power $a \ge 2$.  Note that the sum of the digits in the inside boxes of the horizontal rows add up to $24$.  It is noted here that all shell function calculations that end in the digit $3$ are found in the rows of shell functions with composite powers which yield no prime numbers.\\

\begin{center}
\includegraphics[scale=1]{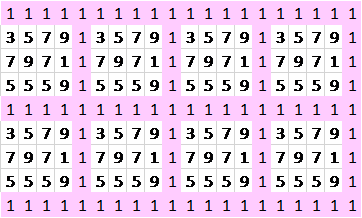}
\end{center}
\begin{center}
\textit{Fig. 6 - Repeating pattern of the last digits of shell calculations for base $n \ge 2$ and power $a \ge 2$ }  
\end{center}

\noindent A fourth observation made from the charts in Figs. 4 and 5 is that there appear to be some rather large "holes" in the distribution of the prime numbers generated by the prime shell function within the range analyzed.  Also, it was observed that there was a rather long delay in the appearance of the first prime numbers for the prime shell functions $n^{37}-(n-1)^{37}$ and $n^{73}-(n-1)^{73}$.\\

\noindent In an effort to gain some level of understanding of the seemingly random occurrences of prime numbers generated by the prime shell function $n^p-(n-1)^p$, one might naturally turn to an application of the prime counting functions introduced by Gauss and Riemann \cite{First50MillionPrimes}.  We recall Gauss' estimate of the number of primes less than some integer limit $x$ as

$$\pi(x) \approx \displaystyle{\frac{x}{ln(x)}}.$$

\noindent Or, as Gauss' more refined prediction allows,
$$\pi(x) \approx Li(x)$$

\begin{center}
where
\end{center}

$$Li(x) = \int_{2}^{x}\frac{1}{\ln(t)} dt.$$

\noindent It is believed that Riemann improved on the logarithmic integral for predicting the number of prime numbers less than or equal to some positive integer limit $x$ (at least for relatively small $x$) by including in Gauss' estimate the powers of primes as illustrated in the following equation \cite{First50MillionPrimes}:

$$\pi(x) +\frac{1}{2}\pi(\sqrt{x})+\frac{1}{3}\pi(\sqrt[3]{x})+...\approx Li(x)$$

\begin{center}
\noindent so that
\end{center}
\emph{•}
$$\pi(x) \approx Li(x)-\frac{1}{2} Li(\sqrt{x})-\frac{1}{3}Li(\sqrt[3]{x})-...=R(x)$$

\noindent in which the last term is the Riemann function more conveniently expressed as

$$R(x)=\sum\limits_{n=1}^\infty \frac{\mu(n)}{n} Li(x^{\frac{1}{n}}).$$

\noindent But the function $Li(x)$ can approximated by the logarithmic sum as follows \cite{First50MillionPrimes}:

$$Li(x) \approx Ls(x) = \frac{1}{\ln(2)}+\frac{1}{\ln(3)}+\frac{1}{\ln(4)}+\frac{1}{\ln(5)}+...+\frac{1}{\ln(x)}.$$

\noindent Calculations that employ the logarithmic sum illustrated above conveniently facilitate the spreadsheet computation and plotting of logarithmic curves that manifest the trend of prime number distribution along the number line of positive integers.  It will be shown that this is no less the case for the calculation of the frequency distribution of the shell prime numbers which, when generated by the prime shell function with the appropriate mathematical operations applied, will approximate the number of shell primes less than or equal to some upper limit $x$ of the base $n$ for some fixed prime power $p$.\\

\noindent In an endeavor to formulate a function parallel to $Li(x)$ which will approximate the frequency and distribution of the shell prime numbers over the realm of the positive integers, it is expedient to introduce a function $M(x)$ that applies both logarithmic and power operations on the prime shell function prior to its integration or summation.  For the integration, we define

$$Mi(x) = \int_{2}^{x}\frac{1}{g[f(n)]}dn$$

\noindent where the function $g$ represents logarithmic and power operations performed on the prime shell function $f(n)$ prior to its integration.  As it turns out, the integration of $g$ of $f(n)$, defined in terms of $f(n)$, will estimate the number of shell primes less than or equal to some upper limit $x$ of the base $n$ using the prime shell function $n^p-(n-1)^p$.  Thus, the following approximation is proposed for calculating the number of shell prime numbers for base $n \leq x$ by the prime shell function:

$$\Pi(x)\approx Mi(x)$$

\noindent where $\Pi(x)$ represents the actual count of prime numbers that exist for the prime shell function for base $n$ less than or equal to some upper limit $x$, and $Mi(x)$ is the aforementioned integral that parallels $Li(x)$ in the Gauss equation.  It will be shown more clearly that $Mi(x)$ is an integral of the prime shell function $f(n)=n^p-(n-1)^p$ which has had logarithmic and power operations performed on it prior to its integration.\\

\noindent Substituting a logarithmic sum $Ms(x)$ in place of the integral $Mi(x)$ in the above equation, similar to how $Ls(x)$ can be used to estimate $Li(x)$, one obtains

$$\Pi(x)\approx Mi(x)\approx Ms(x)=\frac{1}{g(f(2))}+\frac{1}{g(f(3))}+\frac{1}{g(f(4))}+\frac{1}{g(f(5))}+...+\frac{1}{g(f(x))}$$

\noindent where ${g({f(n)})}$ is, up to this point, some undefined function $g$ of the prime shell function $f(n)=n^p-(n-1)^p$.  As it turns out, a good approximation of $Ms(x)$ using the logarithmic summation in place of the integral $Mi(x)$ is found to be

$$Mi(x)\approx Ms(x)=\frac{1}{\sqrt[m]{\ln{f(2)}}}+\frac{1}{\sqrt[m]{\ln{f(3)}}}+\frac{1}{\sqrt[m]{\ln{f(4)}}}+\frac{1}{\sqrt[m]{\ln{f(5)}}}+...+\frac{1}{\sqrt[m]{\ln{f(x)}}}$$

\noindent where the $m^{th}$ root has been empirically determined by the author to be approximately $1.68723$ for the range of base $n \le x$ in this study of the prime shell function $n^p-(n-1)^p$.  Please see Fig. 7 for a graph of the approximating curves generated by the $Ms(x)$ function which were used to estimate the number of shell primes for base $n \le 100$ (bottom curve) and base $n \le 200$ (top curve) compared to the actual count of shell prime numbers $\Pi_{n \le 100}$ and $\Pi_{n \le 200}$, respectively.\\

\begin{center}
\includegraphics[scale=0.55]{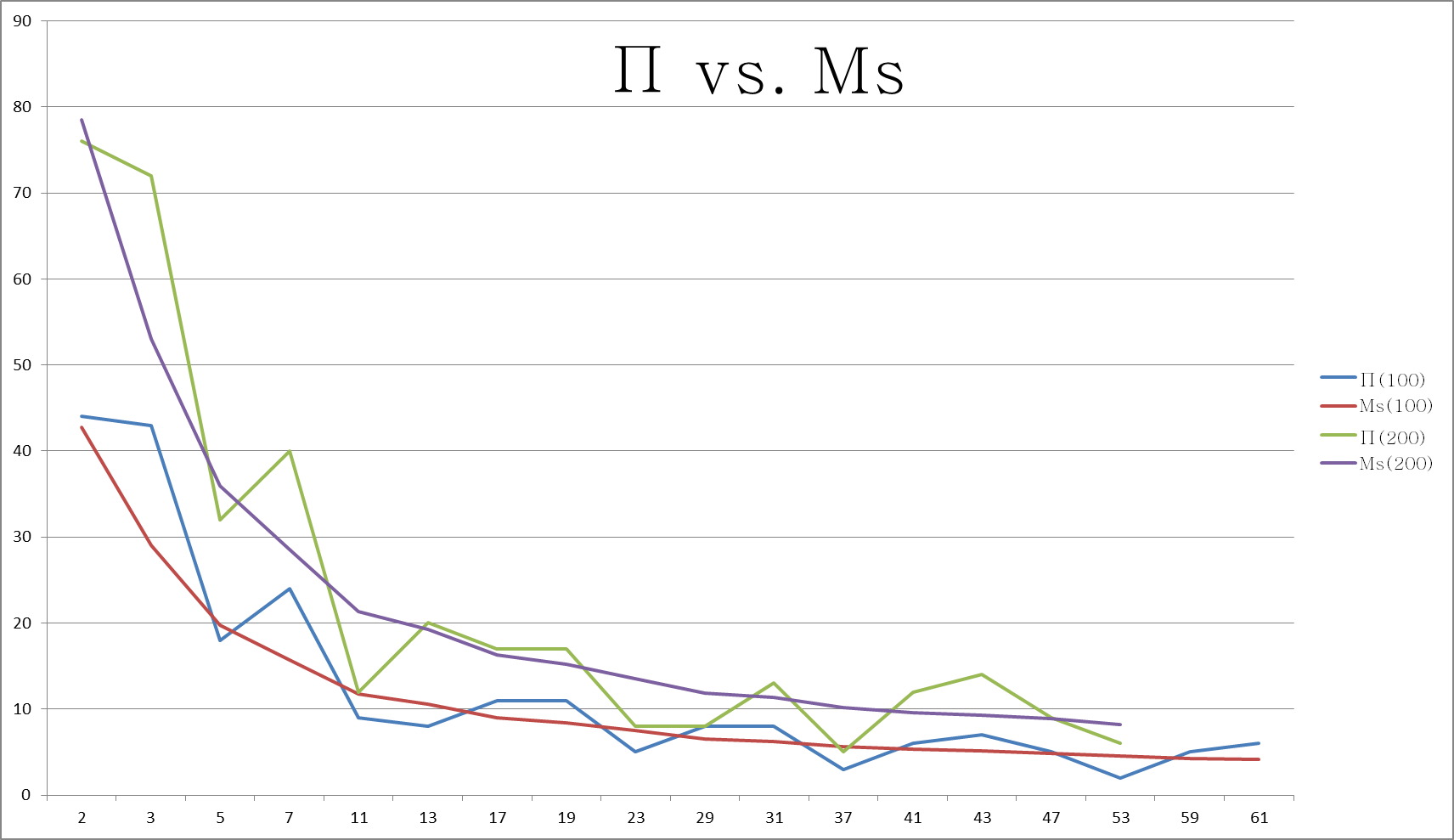}\\
\end{center}

\begin{center}
\textit{Fig. 7 - Graph of the approximating curves generated by the $Ms(x)$ function which estimate the number of shell prime numbers for base $n \le 100$ (bottom curve) and base $n \le 200$ (top curve) for the prime shell function compared to the actual count of shell primes numbers $\Pi_{n \le 100}$ and $\Pi_{n \le 200}$, respectively} 
\end{center}

\noindent It is noted at this point that it would certainly be possible to find constants A and B which could be applied to the $f(n)$ terms of the logarithmic summation function $Ms(x)$ that would appear to provide a better fit of the $Ms(x)$ sum curve to the plot of the actual count of prime numbers within the range of this analysis.  For example, an equation such as $$\frac{1}{A*\sqrt{ln(f(n))}+B}$$ could be formulated to generate a curve that may more closely follow the path of the plot of the actual number of primes generated by the prime shell function within the range of this analysis.  However, it was determined to be more effective to optimize a logarithmic root function to apply to the $f(n)$ terms and omit any additive constants in approximating the average path of the plot of the actual number of primes generated by the prime shell functions.  Although the $f(n)$ terms in the $Ms(x)$ summation that apply constants A and B could create an equation that would appear to be better formed to the points of the actual count of prime numbers less than or equal to some upper limit $x$, $f(n)$ terms that omit any such additive constants but which employ an optimized root of the natural logarithm of the function $f(n)$ will ensure that the curve generated by the $Ms(x)$ sum will remain asymptotic to the x-axis at infinity.  As an example, it was found that an additive constant B could in some cases cause the curve of the $Ms(x)$ function to dip below the x-axis into negative values for the number of estimated primes thus rendering the $Ms(x)$ curve estimate to be non-asymptotic and not a good representation of the diminishing frequency of prime numbers generated as the prime power $p$ of the prime shell function approaches infinity.\\

\noindent Please see Table 1 for the actual count of shell primes $\Pi$ for base $n \le 100$ and base $n \le 200$ for prime powers $p \le 61$ compared to the $Ms(x)$ estimates that employ the $mth$ root of the natural log of its terms.\\

\begin{center}
\begin{tabular}{|c||c|c||c|c|}
\hline
prime power $p$ & $\Pi \le f(100)$ & $Ms \le f(100)$ & $\Pi \le f(200)$ & $Ms \le f(200)$ \\
\hline
2 & 44 & 42.75969 & 76 & 78.48273 \\
\hline
3 & 43 & 29.01307 & 72 & 53.06455 \\
\hline
5 & 18 & 19.71488 & 32 & 35.92022 \\
\hline
7 & 24 & 15.71077 & 40 & 28.56513 \\
\hline
11 & 9 & 11.77596 & 12 & 21.36330 \\
\hline
13 & 8 & 10.61689 & 20 & 19.24779 \\
\hline
17 & 11 & 9.00845 & 17 & 16.31698 \\
\hline
19 & 11 & 8.42015 & 17 & 15.24648 \\
\hline
23 & 5 & 7.50235 & 8 & 13.57796 \\
\hline
29 & 8 & 6.52712 & 8 & 11.80714 \\
\hline
31 & 8 & 6.27155 & 13 & 11.34343 \\
\hline
37 & 3 & 5.64216 & 5 & 10.20206 \\
\hline
41 & 6 & 5.30697 & 12 & 9.59457 \\
\hline
43 & 7 & 5.15842 & 14 & 9.32540 \\
\hline
47 & 5 & 4.89221 & 9 & 8.84316 \\
\hline
53 & 2 & 4.55460 & 6 & 8.23180 \\
\hline
59 & 5 & 4.27319 & --- & --- \\
\hline
61 & 6 & 4.18934 & --- & --- \\
\hline
\end{tabular}
\end{center}

\begin{center}
\textit{Table 1 - Actual count of the shell prime numbers $\Pi$ for base $n \le 100$ and for $\le 200$ for prime powers $p \le 61$ compared to the estimates yielded by the logarithmic sum curves generated with the root value $m = 1.68723$} 
\end{center}

\noindent It was proposed that a good approximation for the $m^{th}$ root of the natural log of the prime shell function $f(n)$ which will produce $Ms(x)$ curves that estimate the number of prime numbers was empirically determined to be 1.68723, and this was illustrated graphically in Fig. 7 for the prime shell function $n^p-(n-1)^p$ for the cases of $x=100$ and $x=200$ for $p \le 61$ and $p \le 53$, respectively.  The author arrived at this empirical estimate by graphing several curves with incremental values of the root $m$ as is illustrated in Fig. 8.  From this graph, one can also determine values of the root $m$ of the natural log of the prime shell function which will generate curves that provide upper and lower bounds for the actual number of primes generated by the prime shell function.  In Fig. 8, the upper bounding curve has a value of $m_{upper}=2.00000$ in the root function $g(f(x))$, and the lower bounding curve has a value of $m_{lower}=1.35759$.  The approximation of $m \approx 1.68723$ was determined to be the best fit $Ms(x)$ curve for the case of $n\le100$, and it is anticipated that this value of $m$ will remain fairly constant as the prime power continues to increase in the prime shell function $n^p-(n-1)^p$ and as the range $x$ of base $n$ in the prime shell function continues toward $\infty$.\\

\begin{center}
\includegraphics[scale=0.73]{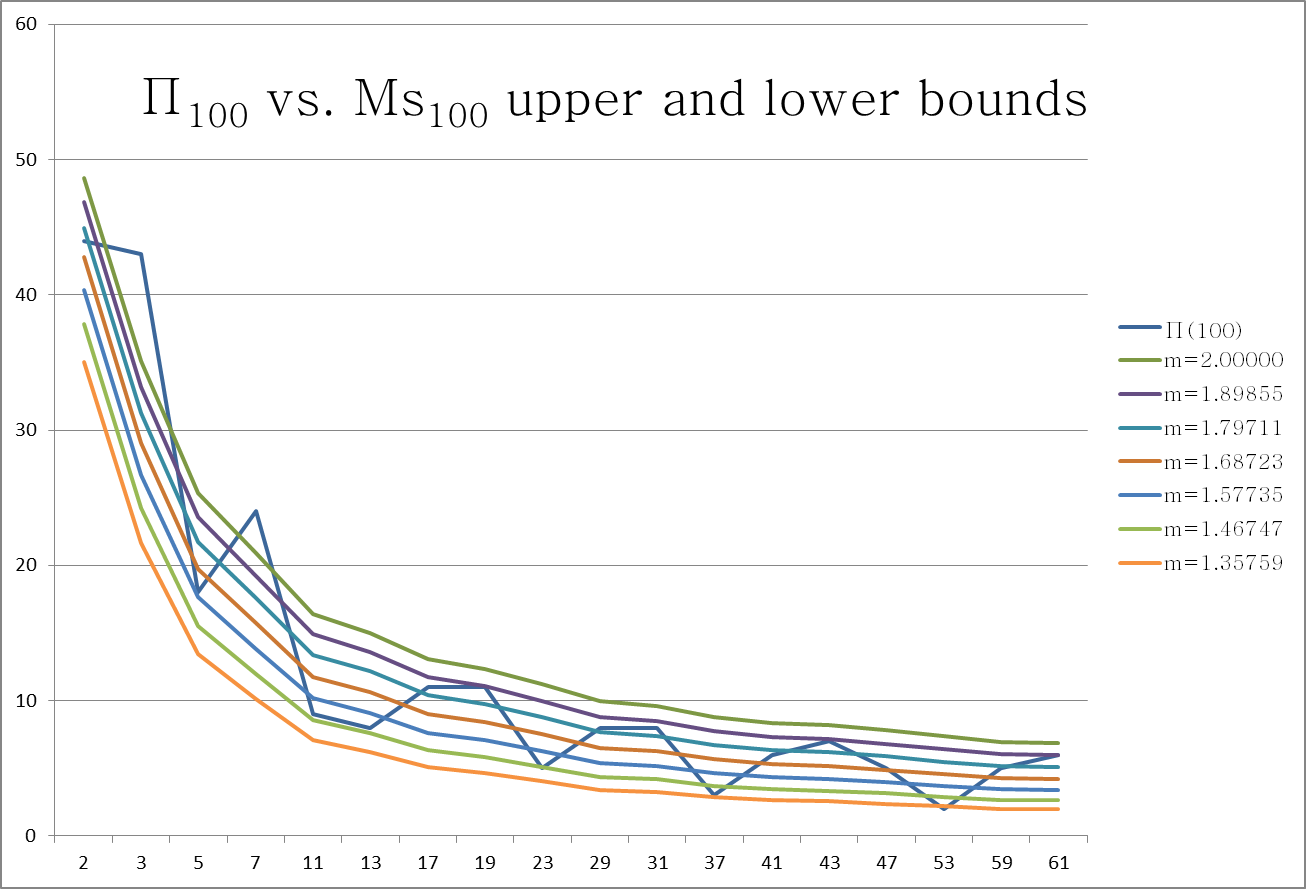}\\
\end{center}

\begin{center}
\textit{Fig. 8 - Logarithmic curves generated by varying the value of the root $m$ to create upper and lower bounds for the estimates of the number of shell primes for base $n \le 100$ in the prime shell function $n^p-(n-1)^p$ for prime power $p \le 61$} 
\end{center}

\noindent If one were to include the powers of primes in the equation for the calculation of the shell primes, similar to what Riemann did to improve the Gauss' $Li$ equation (at least for relatively small $x$), we would get an equation that parallels the Riemann prime counting function such as follows\cite{First50MillionPrimes}:

$$Mi(x) \approx \Pi(x) +\frac{1}{2}\Pi(\sqrt{x})+\frac{1}{3}\Pi(\sqrt[3]{x})+\frac{1}{5}\Pi(\sqrt[5]{x})-\frac{1}{6}\Pi(\sqrt[6]{x})...$$

\noindent which implies

$$\Pi(x)\approx Mi(x) -\frac{1}{2}Mi(\sqrt{x})-\frac{1}{3}Mi(\sqrt[3]{x})-\frac{1}{5}Mi(\sqrt[5]{x})+\frac{1}{6}Mi(\sqrt[6]{x})....$$

\noindent To facilitate the calculation of this improved estimate, one could replace the first few integrals in the above equation with $Ms(x)$ logarithmic sums to get a more accurate approximation of the Riemann method for calculating $\Pi(x)$.  That exercise, however, will be deferred in this research.\\ 

\noindent One who is familiar with Euler's zeta function might naturally ponder the application of that special function to the case of the prime shell function.  We recall the famous equation \cite{HowEulerDiscoveredthezetafunction}

\begin{equation}
\zeta(s) = 1 +\frac{1}{2^s}+\frac{1}{3^s}+\frac{1}{4^s}+\frac{1}{5^s}+... = \frac{1}{1-\frac{1}{2^s}}\cdot\frac{1}{1-\frac{1}{3^s}}\cdot\frac{1}{1-\frac{1}{5^s}}\cdot\frac{1}{1-\frac{1}{7^s}}\cdot...
\end{equation}

\noindent in which Euler sieved the prime numbers out from the sum on the left hand side of Eq. 1 to include in the product on the right-hand side of that equation.  However, when one replaces the integer values in the Euler zeta function in Eq. 1 with the values calculated by the prime shell function $f(n)=n^p-(n-1)^p$, one will find that there is not an efficient way to sieve prime numbers out from the sum on the left hand side of Eq. 1 to create the prime product term on the right-hand side of the equation by either scaling a previous prime number calculation or by employing values of the base $n$ which would predictably make subsequent outputs prime.  So the best one may be able to hope for in this case is to determine some value to be added or subtracted from unity in the numerator on the right-hand side of Eq. 1 to balance both sides of the equation.  Thus, the author introduces the following equation for the non-sieving application of the Euler zeta function to the special case of the prime shell function:

\newcommand{\opM}{\mathop{\vphantom{\sum}\mathchoice
  {\vcenter{\hbox{\huge M}}}
  {\vcenter{\hbox{\Large A}}}{\mathrm{A}}{\mathrm{A}}}\displaylimits}
  
\begin{equation}
\ensuremath \raisebox{-1.5pt}{\scalebox{1.5} {Z(s)}}= \displaystyle\sum\limits_{n=1}^\infty \frac{1}{\Big[n^p-(n-1)^p\Big]^s} = \frac{1+\displaystyle\opM\limits_{i=2}^\infty \Big[(-1)^{i-1}\Big]\displaystyle\sum{}^i}{\displaystyle\prod\limits_{n=2}^\infty {1-\frac{1}{\displaystyle\Big[n^p-(n-1)^p\Big]^s}•}}
\end{equation}

\noindent where

$$\opM\limits_{i=2}^\infty \Big[(-1)^{i-1}\Big]\displaystyle\sum{}^i=-\Sigma_1 \Sigma_2+\Sigma_1 \Sigma_2 \Sigma_3-\Sigma_1 \Sigma_2 \Sigma_3 \Sigma_4+\Sigma_1 \Sigma_2 \Sigma_3 \Sigma_4 \Sigma_5-\Sigma_1 \Sigma_2 \Sigma_3 \Sigma_4 \Sigma_5 \Sigma_6+... \implies$$ 

$$\Sigma_1 \Sigma_2 = \sum\limits_{i=2}^\infty \sum\limits_{j=i}^\infty \frac{1}{f(i)} \frac{1}{f(j)}$$

$$\Sigma_1 \Sigma_2 \Sigma_3 = \sum\limits_{i=2}^\infty \sum\limits_{j=i}^\infty \sum\limits_{k=j+1}^\infty \frac{1}{f(i)} \frac{1}{f(j)} \frac{1}{f(k)}$$

$$\Sigma_1 \Sigma_2 \Sigma_3 \Sigma_4 = \sum\limits_{i=2}^\infty \sum\limits_{j=i}^\infty \sum\limits_{k=j+1}^\infty \sum\limits_{l=k+1}^\infty \frac{1}{f(i)} \frac{1}{f(j)} \frac{1}{f(k)} \frac{1}{f(l)}$$

$$\Sigma_1 \Sigma_2 \Sigma_3 \Sigma_4 \Sigma_5 = \sum\limits_{i=2}^\infty \sum\limits_{j=i}^\infty \sum\limits_{k=j+1}^\infty \sum\limits_{l=k+1}^\infty \sum\limits_{m=l+1}^\infty \frac{1}{f(i)} \frac{1}{f(j)} \frac{1}{f(k)} \frac{1}{f(l)}\frac{1}{f(m)}$$

$$\Sigma_1 \Sigma_2 \Sigma_3 \Sigma_4 \Sigma_5 \Sigma_6= \sum\limits_{i=2}^\infty \sum\limits_{j=i}^\infty \sum\limits_{k=j+1}^\infty \sum\limits_{l=k+1}^\infty \sum\limits_{m=l+1}^\infty \sum\limits_{n=m+1}^\infty \frac{1}{f(i)} \frac{1}{f(j)} \frac{1}{f(k)} \frac{1}{f(l)} \frac{1}{f(m)} \frac{1}{f(n)}$$

\begin{center}
\bf.
\end{center}
\begin{center}
\bf.
\end{center}
\begin{center}
\bf.
\end{center}

\noindent As intimated by the first row of Fig. 9, if one applies the realm of positive integers $i$ to the Euler zeta function using the non-sieving  method so that the product term in the denominator on the right-hand side of Eq. 1 includes all the prime and composite terms, then the numerator on the right-hand side of Eq. 1 will tend to zero as the limit $x$ tends to infinity.  Thus, for the case of the positive integers, one might associate the value of all the composite terms in the product term in the denominator of the right-hand side of the non-sieving Euler zeta function, if they were included, to $-1$ in the numerator on that side of the equation, because when all the composite terms are included in the product term in the denominator of Eq. 1, the numerator on the right-hand side of that equation tends to zero (i.e., $1+(M)=0$) as $x$ approaches $\infty$.  Compare that to the value of unity attained in the numerator of the right-hand side Eq. 1 when all the composite terms are sieved out of the denominator of the product term on the right-hand side of the equation as Euler did when he discovered this famous relationship between the composite and prime numbers.  So one might postulate that the total effect of sieving the composite terms out of the product term in the denominator of the right-hand side of the Euler zeta function in Eq. 1 is that the numerator on the right-hand side of the equation changes from zero to unity as the composite-numbered terms are eliminated.  It is observed that while the "M-series"

$$\displaystyle\opM\limits_{i=2}^\infty \Big{[(-1)^{i-1}\Big]\displaystyle\sum{}^i}$$

\noindent in the numerator of the right-hand side of the non-sieving zeta function represents an infinite sum of infinite sums, its value is constrained between $0$ and $-1$ throughout the application of this function to the realm of positive integers as well as to the range of integer calculations yielded by the prime shell function $n^p-(n-1)^p$.

\begin{center}
\includegraphics[scale=0.6]{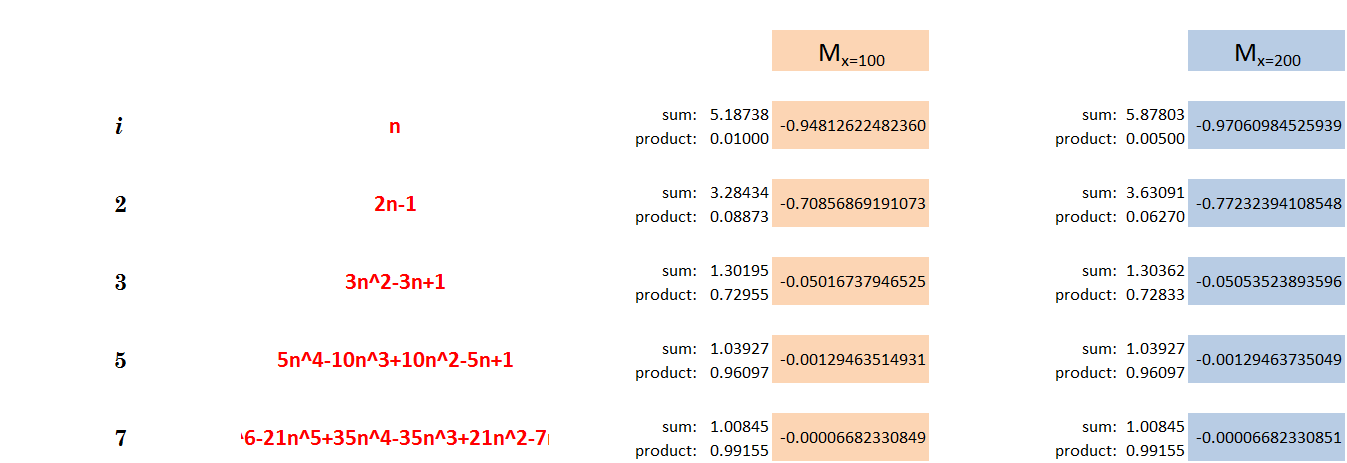}
\end{center}

\begin{center}
\textit{Fig. 9 - M-series calculations for the realm of the positive integers $i$ and for the first four prime powers $p=2,3,5,7$ of the prime shell function $n^p-(n-1)^p$ for $x=100$ and for $x=200$} 
\end{center}

\noindent In Figure 9, descending rows below row $i$ depict the M-series values $M_{x=100}$ and $M_{x=200}$, respectively, that result from the non-sieving application of the zeta function to the prime shell function of powers $2, 3, 5$ and $7$, respectively.  Included to the left of each M-series value are the $Ms$ (sum) and $\Pi$ (product) calculations for that row of the prime shell function that were used to calculate the M-series values under the columns $M_{x=100}$ and $M_{x=200}$.  Fig. 10 is a graph plotted that superimposes the trends of $M_{x=100}$ and $M_{x=200}$ for the case of the non-sieving application of the Euler zeta function, and it illustrates the change in the value of the M-series as the prime power in the prime shell function increases from $2$ to $7$ and as the range of the base $n$ of the prime shell function increases from $x=100$ to $x=200$.  Notice the trend is that the M-series increases from $-1$ to $0$ asymptotically as the power of the prime shell equation increases and that the increase in the M-series value between $x=100$ to $x=200$ is very slight.  It is expected that as $x$ approaches $\infty$ for row $i$ in the case of the application of the non-sieving zeta function to the realm of positive integers that the value of the M-series will approach $-1$ thereby rendering the numerator of the right-hand side of Eq. 2 to zero.  It is also expected that the limit of the M-series in the descending rows below row $i$ which summarize the results of the non-sieving application of the Euler zeta function to the prime shell functions of powers $2, 3, 5$ and $7$ will approach stable limits as the range $x$ of the base $n$ approaches $\infty$.  From the trends manifested by the graphs in Fig. 10, it is anticipated that the M-series values will approach the following approximate limits as $x$ approaches $\infty$:\\

$$M_{{i}_{x=\infty}}\approx -1$$ $$M_{{2}_{(x=\infty)}}\approx -0.8$$ $$M_{{3}_{(x=\infty)}}\approx -0.051$$ $$M_{{5}_{(x=\infty)}}\approx -0.0013$$ $$M_{{7}_{(x=\infty)}}\approx -0.000067$$

\begin{center}
\includegraphics[scale=1.1]{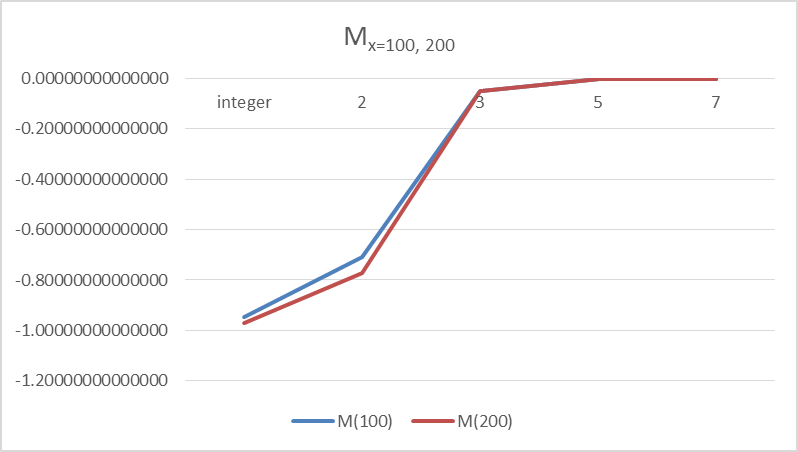}
\end{center}
\begin{center}
\textit{Fig. 10 - Graph of the M-series superimposed for $x = 100$ and $x = 200$ for \textit{i}$=$integer and for the first four prime powers $p=2,3,5,7$ of the shell prime function.  Note that the M-series is asymptotic to $0$ as the prime power $p$ in the prime shell function $n^p-(n-1)^p$ increases from $2$ to $7$.} \\
\end{center}

\noindent Table 2 tabulates the values of $M_{x=100}$ and $M_{x=200}$ for the non-sieving application of the Euler zeta function to the positive integers $i$ and to the prime shell function with powers $2, 3, 5$ and $7$.  It is anticipated that such calculations will be able to tell us something about the density of shell prime numbers less than a given limit $x$ of the range of base $n$ for any given prime power $p$ in the prime shell function.  Further study should reveal more insight into the properties of this special series and its potential for predicting the number of primes generated by the prime shell function and for any other integer-valued polynomial in general which has the capability to generate prime numbers and which has as its input domain the realm of positive integers.

\begin{center}
\begin{tabular}{|c||c||c|}
\hline
{} & ${\displaystyle\opM\limits_{}^{}}_{x=100}$ & ${\displaystyle\opM\limits_{}^{}}_{x=200}$ \\
\hline
integer $i$ & -0.94812622482360
 & -0.97060984525939
 \\
\hline
$p$ = $2$ & -0.70856869191073
 & -0.77232394108548
 \\
\hline
$p$ = $3$ & -0.05016737946525
 & -0.05053523893596
 \\
\hline
$p$ = $5$ & -0.00129463514931
 & -0.00129463735049
 \\
\hline
$p$ = $7$ & -0.00006682330849
 & -0.00006682330851
 \\
\hline
\end{tabular}
\end{center}

\begin{center}
\textit{Table 2} \\
\end{center}

\noindent As an exercise to promote the notion of a non-sieving application of the Euler zeta function to integer-valued polynomial functions in general, we shall apply the non-sieving method to any integer-generating polynomial for which the input domain of the function is the positive integers that generate successive integer values by the polynomial for processing in the zeta function.  Thus, we introduce the following "big Zeta" equation for this generalized application: 

\begin{equation}
\ensuremath \raisebox{-1.5pt}{\scalebox{1.5} {Z}} \displaystyle = \displaystyle {1+} \displaystyle\sum\limits_{n=1}^\infty \frac{1}{f(n)} = \frac{1+\displaystyle\opM\limits_{i=1}^\infty (-1)^{i}\displaystyle\sum{}^{i+1}}{\displaystyle\prod\limits_{n=1}^\infty {1-\frac{1}{\displaystyle f(n)}•}}.
\end{equation}

\noindent where $f(n)$ represents any polynomial function which yields successively increasing positive integer values when the input domain of the function is the positive integers and where $f(1) \neq 1$.\\

\noindent We begin our example of the non-sieving application of the Euler zeta function to integer-valued polynomials in general by writing

\begin{equation}
Z = 1 +\frac{1}{f(1)}+\frac{1}{f(2)}+\frac{1}{f(3)}+\frac{1}{f(4)}+\frac{1}{f(5)}+\frac{1}{f(6)}+\frac{1}{f(7)}+\frac{1}{f(8)}+\frac{1}{f(9)}+... .
\end{equation}

\noindent We first multiply both sides of Eq. 4 by the first fraction on the right hand side of that equation which yields

\begin{multline*}
$$\frac{1}{f(1)}\cdot Z = \frac{1}{f(1)}+\frac{1}{f(1)}\cdot\frac{1}{f(1)}+\frac{1}{f(1)}\cdot\frac{1}{f(2)}+\frac{1}{f(1)}\cdot\frac{1}{f(3)}+\frac{1}{f(1)}\cdot\frac{1}{f(4)}+
\frac{1}{f(1)}\cdot\frac{1}{f(5)}+\frac{1}{f(1)}\cdot\frac{1}{f(6)}+\frac{1}{f(1)}\cdot\frac{1}{f(7)}+\frac{1}{f(1)}\cdot\frac{1}{f(8)}+\frac{1}{f(1)}\cdot\frac{1}{f(9)}+...$$
\end{multline*}

\noindent and when we subtract this result from both sides of Eq. 4, we get

\begin{equation}
\Big(1-\frac{1}{f(1)}\Big)\cdot Z =1 +\frac{1}{f(2)}+\frac{1}{f(3)}+\frac{1}{f(4)}+\frac{1}{f(5)}+\frac{1}{f(6)}+\frac{1}{f(7)}+\frac{1}{f(8)}+\frac{1}{f(9)}+...-\displaystyle\sum\limits_{n=1}^\infty \frac{1}{f(1) \cdot f(n)}.
\end{equation}

\noindent Continuing with the multiplication of Eq. 5 by the next fraction in the sequence, $\frac{1}{f(2)}$: 

\begin{multline*}
$$\frac{1}{f(2)}\cdot\Big(1-\frac{1}{f(1)}\Big)\cdot Z =\frac{1}{f(2)}+\frac{1}{f(2)}\cdot\frac{1}{f(2)}+\frac{1}{f(2)}\cdot\frac{1}{f(3)}+\frac{1}{f(2)}\cdot\frac{1}{f(4)}+\frac{1}{f(2)}\cdot\frac{1}{f(5)}+\frac{1}{f(2)}\cdot\frac{1}{f(6)}+\frac{1}{f(2)}\cdot\frac{1}{f(7)}+\frac{1}{f(2)}\cdot\frac{1}{f(8)}+\frac{1}{f(2)}\cdot\frac{1}{f(9)}+...\\-\displaystyle\sum\limits_{n=1}^\infty \frac{1}{f(1) \cdot f(2)\cdot f(n)}$$
\end{multline*}

\noindent so that when we subtract this result from both sides of Eq. 5, we end up with

\begin{multline*}
$$\Big(1-\frac{1}{f(2)}\Big)\Big(1-\frac{1}{f(1)}\Big)\cdot Z =1 +\frac{1}{f(3)}+\frac{1}{f(4)}+\frac{1}{f(5)}+\frac{1}{f(6)}+\frac{1}{f(7)}+\frac{1}{f(8)}+\frac{1}{f(9)}+...\\-\displaystyle\sum\limits_{n=1}^\infty \frac{1}{f(1) \cdot f(n)}-\displaystyle\sum\limits_{n=2}^\infty \frac{1}{f(2) \cdot f(n)}+\displaystyle\sum\limits_{n=1}^\infty \frac{1}{f(1) \cdot f(2) \cdot f(n)}.$$
\end{multline*}

\noindent Continuing with the multiplication of this new equation by the next term in the sequence, $\frac{1}{f(3)}$:

\begin{multline*}
$$\frac{1}{f(3)}\cdot\Big(1-\frac{1}{f(2)}\Big)\cdot\Big(1-\frac{1}{f(1)}\Big)\cdot Z =\frac{1}{f(3)}+\frac{1}{f(3)}\cdot\frac{1}{f(3)}+\frac{1}{f(3)}\cdot\frac{1}{f(4)}+\frac{1}{f(3)}\cdot\frac{1}{f(5)}+\frac{1}{f(3)}\cdot\frac{1}{f(6)}+\frac{1}{f(3)}\cdot\frac{1}{f(7)}+\frac{1}{f(3)}\cdot\frac{1}{f(8)}+\frac{1}{f(3)}\cdot\frac{1}{f(9)}+...\\-\displaystyle\sum\limits_{n=1}^\infty \frac{1}{f(1) \cdot f(3)\cdot f(n)}-\displaystyle\sum\limits_{n=2}^\infty \frac{1}{f(2) \cdot f(3)\cdot f(n)}+\displaystyle\sum\limits_{n=1}^\infty \frac{1}{f(1) \cdot f(2)\cdot f(3)\cdot f(n)}$$
\end{multline*}

\noindent so that

\begin{multline*}
$$\Big(1-\frac{1}{f(3)}\Big)\Big(1-\frac{1}{f(2)}\Big)\Big(1-\frac{1}{f(1)}\Big)\cdot Z =1 +\frac{1}{f(4)}+\frac{1}{f(5)}+\frac{1}{f(6)}+\frac{1}{f(7)}+\frac{1}{f(8)}+\frac{1}{f(9)}+...\\-\displaystyle\sum\limits_{n=1}^\infty \frac{1}{f(1) \cdot f(n)}-\displaystyle\sum\limits_{n=2}^\infty \frac{1}{f(2) \cdot f(n)}-\displaystyle\sum\limits_{n=3}^\infty \frac{1}{f(3) \cdot f(n)}\\+\displaystyle\sum\limits_{n=1}^\infty \frac{1}{f(1)\cdot f(2) \cdot f(n)}+\displaystyle\sum\limits_{n=1}^\infty \frac{1}{f(1) \cdot f(3)\cdot f(n)}+\displaystyle\sum\limits_{n=2}^\infty \frac{1}{f(2) \cdot f(3)\cdot f(n)}-\displaystyle\sum\limits_{n=1}^\infty \frac{1}{f(1) \cdot f(2)\cdot f(3)\cdot f(n)}.$$
\end{multline*}

\noindent Continuing with the multiplication of this new equation by the next term in the sequence, $\frac{1}{f(4)}$:

\begin{multline*}
$$\frac{1}{f(4)}\cdot\Big(1-\frac{1}{f(3)}\Big)\cdot\Big(1-\frac{1}{f(2)}\Big)\cdot\Big(1-\frac{1}{f(1)}\Big)\cdot Z =\frac{1}{f(4)}+\frac{1}{f(4)}\cdot\frac{1}{f(4)}+\frac{1}{f(4)}\cdot\frac{1}{f(5)}+\frac{1}{f(4)}\cdot\frac{1}{f(6)}+\frac{1}{f(4)}\cdot\frac{1}{f(7)}+\frac{1}{f(4)}\cdot\frac{1}{f(8)}+\frac{1}{f(4)}\cdot\frac{1}{f(9)}+...\\-\displaystyle\sum\limits_{n=1}^\infty \frac{1}{f(1)\cdot f(4) \cdot f(n)}-\displaystyle\sum\limits_{n=2}^\infty \frac{1}{f(2) \cdot f(4)\cdot f(n)}-\displaystyle\sum\limits_{n=3}^\infty \frac{1}{f(3) \cdot f(4)\cdot f(n)}+\displaystyle\sum\limits_{n=1}^\infty \frac{1}{f(1)\cdot f(2) \cdot f(4)\cdot f(n)}\\+\displaystyle\sum\limits_{n=1}^\infty \frac{1}{f(1) \cdot f(3)\cdot f(4)\cdot f(n)}+\displaystyle\sum\limits_{n=2}^\infty \frac{1}{f(2) \cdot f(3)\cdot f(4)\cdot f(n)}-\displaystyle\sum\limits_{n=1}^\infty \frac{1}{f(1) \cdot f(2)\cdot f(3)\cdot f(4)\cdot f(n)}$$
\end{multline*}

\noindent so that

\begin{multline*}
$$\Big(1-\frac{1}{f(4)}\Big)\Big(1-\frac{1}{f(3)}\Big)\Big(1-\frac{1}{f(2)}\Big)\Big(1-\frac{1}{f(1)}\Big)\cdot Z =1 +\frac{1}{f(5)}+\frac{1}{f(6)}+\frac{1}{f(7)}+\frac{1}{f(8)}+\frac{1}{f(9)}+...\\-\displaystyle\sum\limits_{n=1}^\infty \frac{1}{f(1) \cdot f(n)}-\displaystyle\sum\limits_{n=2}^\infty \frac{1}{f(2) \cdot f(n)}-\displaystyle\sum\limits_{n=3}^\infty \frac{1}{f(3) \cdot f(n)}-\displaystyle\sum\limits_{n=4}^\infty \frac{1}{f(4) \cdot f(n)}\\+\displaystyle\sum\limits_{n=1}^\infty \frac{1}{f(1)\cdot f(2) \cdot f(n)}+\displaystyle\sum\limits_{n=1}^\infty \frac{1}{f(1) \cdot f(3)\cdot f(n)}+\displaystyle\sum\limits_{n=2}^\infty \frac{1}{f(2) \cdot f(3)\cdot f(n)}-\displaystyle\sum\limits_{n=1}^\infty \frac{1}{f(1) \cdot f(2)\cdot f(3)\cdot f(n)}\\+\displaystyle\sum\limits_{n=1}^\infty \frac{1}{f(1)\cdot f(4) \cdot f(n)}+\displaystyle\sum\limits_{n=2}^\infty \frac{1}{f(2) \cdot f(4)\cdot f(n)}+\displaystyle\sum\limits_{n=3}^\infty \frac{1}{f(3) \cdot f(4)\cdot f(n)}-\displaystyle\sum\limits_{n=1}^\infty \frac{1}{f(1)\cdot f(2) \cdot f(4)\cdot f(n)}\\-\displaystyle\sum\limits_{n=1}^\infty \frac{1}{f(1) \cdot f(3)\cdot f(4)\cdot f(n)}-\displaystyle\sum\limits_{n=2}^\infty \frac{1}{f(2) \cdot f(3)\cdot f(4)\cdot f(n)}+\displaystyle\sum\limits_{n=1}^\infty \frac{1}{f(1) \cdot f(2)\cdot f(3)\cdot f(4)\cdot f(n)}.$$
\end{multline*}

\noindent Continuing with the multiplication of of this new equation by the next term in the sequence, $\frac{1}{f(5)}$:

\begin{multline*}
$$\frac{1}{f(5)}\cdot\Big(1-\frac{1}{f(4)}\Big)\cdot\Big(1-\frac{1}{f(3)}\Big)\cdot\Big(1-\frac{1}{f(2)}\Big)\cdot\Big(1-\frac{1}{f(1)}\Big)\cdot Z =\frac{1}{f(5)}+\frac{1}{f(5)}\cdot\frac{1}{f(5)}+\frac{1}{f(5)}\cdot\frac{1}{f(6)}+\frac{1}{f(5)}\cdot\frac{1}{f(7)}+\frac{1}{f(5)}\cdot \frac{1}{f(8)}+\frac{1}{f(5)}\cdot \frac{1}{f(9)}+...\\-\displaystyle\sum\limits_{n=1}^\infty \frac{1}{f(1) \cdot f(5) \cdot f(n)}-\displaystyle\sum\limits_{n=2}^\infty \frac{1}{f(2) \cdot f(5) \cdot f(n)}-\displaystyle\sum\limits_{n=3}^\infty \frac{1}{f(3) \cdot f(5)\cdot f(n)}\\-\displaystyle\sum\limits_{n=4}^\infty \frac{1}{f(4) \cdot f(5) \cdot f(n)}+\displaystyle\sum\limits_{n=1}^\infty \frac{1}{f(1)\cdot f(2) \cdot f(5)\cdot f(n)}+\displaystyle\sum\limits_{n=1}^\infty \frac{1}{f(1) \cdot f(3)\cdot f(5)\cdot f(n)}\\+\displaystyle\sum\limits_{n=2}^\infty \frac{1}{f(2) \cdot f(3)\cdot f(5) \cdot f(n)}-\displaystyle\sum\limits_{n=1}^\infty \frac{1}{f(1) \cdot f(2)\cdot f(3) \cdot f(5) \cdot f(n)}+\displaystyle\sum\limits_{n=1}^\infty \frac{1}{f(1) \cdot f(4) \cdot f(5) \cdot f(n)}\\+\displaystyle\sum\limits_{n=2}^\infty \frac{1}{f(2) \cdot f(4) \cdot f(5) \cdot f(n)}+\displaystyle\sum\limits_{n=3}^\infty \frac{1}{f(3) \cdot f(4) \cdot f(5) \cdot f(n)}-\displaystyle\sum\limits_{n=1}^\infty \frac{1}{f(1) \cdot f(2) \cdot f(4)\cdot f(5) \cdot f(n)}\\-\displaystyle\sum\limits_{n=1}^\infty \frac{1}{f(1) \cdot f(3) \cdot f(4) \cdot f(5) \cdot f(n)}-\displaystyle\sum\limits_{n=2}^\infty \frac{1}{f(2) \cdot f(3) \cdot f(4) \cdot f(5) \cdot f(n)}+\displaystyle\sum\limits_{n=1}^\infty \frac{1}{f(1) \cdot f(2) \cdot f(3) \cdot f(4) \cdot f(5) \cdot f(n)}$$
\end{multline*}

\noindent so that

\begin{multline*}
$$\Big(1-\frac{1}{f(5)}\Big)\Big(1-\frac{1}{f(4)}\Big)\Big(1-\frac{1}{f(3)}\Big)\Big(1-\frac{1}{f(2)}\Big)\Big(1-\frac{1}{f(1)}\Big)\cdot Z =1 +\frac{1}{f(6)}+\frac{1}{f(7)}+\frac{1}{f(8)}+\frac{1}{f(9)}+...\\-\displaystyle\sum\limits_{n=1}^\infty \frac{1}{f(1) \cdot f(n)}-\displaystyle\sum\limits_{n=2}^\infty \frac{1}{f(2) \cdot f(n)}-\displaystyle\sum\limits_{n=3}^\infty \frac{1}{f(3) \cdot f(n)}-\displaystyle\sum\limits_{n=4}^\infty \frac{1}{f(4) \cdot f(n)}-\displaystyle\sum\limits_{n=5}^\infty \frac{1}{f(5) \cdot f(n)}\\+\displaystyle\sum\limits_{n=1}^\infty \frac{1}{f(1)\cdot f(2) \cdot f(n)}+\displaystyle\sum\limits_{n=1}^\infty \frac{1}{f(1) \cdot f(3)\cdot f(n)}+\displaystyle\sum\limits_{n=2}^\infty \frac{1}{f(2) \cdot f(3)\cdot f(n)}-\displaystyle\sum\limits_{n=1}^\infty \frac{1}{f(1) \cdot f(2)\cdot f(3)\cdot f(n)}\\+\displaystyle\sum\limits_{n=1}^\infty \frac{1}{f(1)\cdot f(4) \cdot f(n)}+\displaystyle\sum\limits_{n=2}^\infty \frac{1}{f(2) \cdot f(4)\cdot f(n)}+\displaystyle\sum\limits_{n=3}^\infty \frac{1}{f(3) \cdot f(4)\cdot f(n)}-\displaystyle\sum\limits_{n=1}^\infty \frac{1}{f(1)\cdot f(2) \cdot f(4)\cdot f(n)}\\-\displaystyle\sum\limits_{n=1}^\infty \frac{1}{f(1) \cdot f(3)\cdot f(4)\cdot f(n)}-\displaystyle\sum\limits_{n=2}^\infty \frac{1}{f(2) \cdot f(3)\cdot f(4)\cdot f(n)}+\displaystyle\sum\limits_{n=1}^\infty \frac{1}{f(1) \cdot f(2)\cdot f(3) \cdot f(4) \cdot f(n)}\\+\displaystyle\sum\limits_{n=1}^\infty \frac{1}{f(1) \cdot f(5) \cdot f(n)}+\displaystyle\sum\limits_{n=2}^\infty \frac{1}{f(2) \cdot f(5) \cdot f(n)}+\displaystyle\sum\limits_{n=3}^\infty \frac{1}{f(3) \cdot f(5)\cdot f(n)}+\displaystyle\sum\limits_{n=4}^\infty \frac{1}{f(4) \cdot f(5) \cdot f(n)}\\-\displaystyle\sum\limits_{n=1}^\infty \frac{1}{f(1)\cdot f(2) \cdot f(5)\cdot f(n)}-\displaystyle\sum\limits_{n=1}^\infty \frac{1}{f(1) \cdot f(3)\cdot f(5)\cdot f(n)}-\displaystyle\sum\limits_{n=2}^\infty \frac{1}{f(2) \cdot f(3)\cdot f(5) \cdot f(n)}\\+\displaystyle\sum\limits_{n=1}^\infty \frac{1}{f(1) \cdot f(2)\cdot f(3) \cdot f(5) \cdot f(n)}-\displaystyle\sum\limits_{n=1}^\infty \frac{1}{f(1) \cdot f(4) \cdot f(5) \cdot f(n)}-\displaystyle\sum\limits_{n=2}^\infty \frac{1}{f(2) \cdot f(4) \cdot f(5) \cdot f(n)}\\-\displaystyle\sum\limits_{n=3}^\infty \frac{1}{f(3) \cdot f(4) \cdot f(5) \cdot f(n)}+\displaystyle\sum\limits_{n=1}^\infty \frac{1}{f(1) \cdot f(2) \cdot f(4)\cdot f(5) \cdot f(n)}+\displaystyle\sum\limits_{n=1}^\infty \frac{1}{f(1) \cdot f(3) \cdot f(4) \cdot f(5) \cdot f(n)}\\+\displaystyle\sum\limits_{n=2}^\infty \frac{1}{f(2) \cdot f(3) \cdot f(4) \cdot f(5) \cdot f(n)}-\displaystyle\sum\limits_{n=1}^\infty \frac{1}{f(1) \cdot f(2) \cdot f(3) \cdot f(4) \cdot f(5) \cdot f(n)}.$$
\end{multline*}\\\\

\noindent We continue in this fashion until all the fraction terms on the right-hand side of the big Zeta function in Eq. 4 are eliminated.  This, of course, will only occur at $\infty$.  But as we continue to eliminate the terms on the right-hand side of the big Zeta function in Eq. 4, a value for what we will call the "M-series" function

$$\opM\limits_{i=1}^\infty (-1)^{i}\displaystyle\sum{}^{i+1}$$

\noindent emerges which adds to unity in the numerator on the right-hand side of Eq. 3.  It is observed that the value of this M-series function increases between the bounds of $-1$ and $0$ as the terms, both composite and prime, are eliminated from the right-hand side of the big Zeta function Eq. 4.  To evaluate this M-series for the operations we have performed thus far with the big Zeta function on the generalized polynomial function $f(n)$ in this example, we rearrange the summation terms to obtain\\

\begin{multline*}
$$\displaystyle\sum\limits_{n=1}^\infty \frac{1}{f(n)}\cdot \Big(\color{red}-\frac{1}{f(1)}\color{black}+\frac{1}{f(1)\cdot f(2)}+\frac{1}{f(1)\cdot f(3)}-\frac{1}{f(1)\cdot f(2)\cdot f(3)}+\frac{1}{f(1)\cdot f(4)}-\frac{1}{f(1)\cdot f(2)\cdot f(4)}-\frac{1}{f(1)\cdot f(3)\cdot f(4)}\\+\frac{1}{f(1)\cdot f(2)\cdot f(3)\cdot f(4)}+\frac{1}{f(1)\cdot f(5)}-\frac{1}{f(1)\cdot f(2)\cdot f(5)}-\frac{1}{f(1)\cdot f(3)\cdot f(5)}+\frac{1}{f(1)\cdot f(2)\cdot f(3)\cdot f(5)}\\-\frac{1}{f(1)\cdot f(4)\cdot f(5)}+\frac{1}{f(1)\cdot f(2)\cdot f(4)\cdot f(5)}+\frac{1}{f(1)\cdot f(3)\cdot f(4)\cdot f(5)}-\frac{1}{f(1)\cdot f(2)\cdot f(3)\cdot f(4)\cdot f(5)}...\Big)$$
\end{multline*}

\begin{multline*}
$$\quad\displaystyle\sum\limits_{n=2}^\infty \frac{1}{f(n)} \cdot\Big(\color{red}-\frac{1}{f(2)}\color{black}+\frac{1}{f(2)\cdot f(3)}+\frac{1}{f(2)\cdot f(4)}-\frac{1}{f(2)\cdot f(3)\cdot f(4)}+\frac{1}{f(2)\cdot f(5)}-\frac{1}{f(2)\cdot f(3)\cdot f(5)}\\-\frac{1}{f(2)\cdot f(4)\cdot f(5)}+\frac{1}{f(2)\cdot f(3)\cdot f(4)\cdot f(5)}...\Big)$$
\end{multline*}

\begin{flalign*}
&\quad \displaystyle\sum\limits_{n=3}^\infty \frac{1}{f(n)} \cdot\Big(\color{red}-\frac{1}{f(3)}\color{black}+\frac{1}{f(3)\cdot f(4)}+\frac{1}{f(3)\cdot f(5)}-\frac{1}{f(3)\cdot f(4)\cdot f(5)}...\Big)&
\end{flalign*}

\begin{flalign*}
&\quad \displaystyle\sum\limits_{n=4}^\infty \frac{1}{f(n)}\cdot \Big(\color{red}-\frac{1}{f(4)}\color{black}+\frac{1}{f(4)\cdot f(5)}...\Big)&
\end{flalign*}

\begin{flalign*}
&\quad \displaystyle\sum\limits_{n=5}^\infty \frac{1}{f(n)}\cdot \Big(\color{red}-\frac{1}{f(5)}\color{black}...\Big)&
\end{flalign*}

\begin{flalign*}
&\quad \displaystyle\sum\limits_{n=6}^\infty \frac{1}{f(n)}.&
\end{flalign*}

\noindent Once the summation terms are thus grouped in this fashion according to their lower limits, we begin collecting like terms across the groups to organize the sums according to the number of terms which their product will contain when they are expanded with their coefficients.  Referring to the group above, we begin by extracting the summations with the fewest terms which have been highlighted in \color{red} red \color{black}, i.e.,

\begin{flalign*}
&\quad \color{red}-\frac{1}{f(1)}\color{black}\cdot\displaystyle\sum\limits_{n=1}^\infty \frac{1}{f(n)}\color{red}-\frac{1}{f(2)}\color{black}\cdot\displaystyle\sum\limits_{n=2}^\infty \frac{1}{f(n)}\color{red}-\frac{1}{f(3)}\color{black}\cdot\displaystyle\sum\limits_{n=3}^\infty \frac{1}{f(n)}\color{red}-\frac{1}{f(4)}\color{black}\cdot\displaystyle\sum\limits_{n=4}^\infty \frac{1}{f(n)}\color{red}-\frac{1}{f(5)}\color{black}\cdot\displaystyle\sum\limits_{n=5}^\infty \frac{1}{f(n)}\color{red}-\color{black}\cdots=&
\end{flalign*}

\begin{flalign*}
&\quad \color{red}-\frac{1}{f(1)}\color{black}\cdot\Big(\frac{1}{f(1)}+\frac{1}{f(2)}+\frac{1}{f(3)}+\frac{1}{f(4)}+\frac{1}{f(5)}+\cdots\Big)&
\end{flalign*}

\begin{flalign*}
&\quad \color{red}-\frac{1}{f(2)}\color{black}\cdot\Big(\frac{1}{f(2)}+\frac{1}{f(3)}+\frac{1}{f(4)}+\frac{1}{f(5)}+\cdots\Big)&
\end{flalign*}

\begin{flalign*}
&\quad \color{red}-\frac{1}{f(3)}\color{black}\cdot\Big(\frac{1}{f(3)}+\frac{1}{f(4)}+\frac{1}{f(5)}+\cdots\Big)&
\end{flalign*}

\begin{flalign*}
&\quad \color{red}-\frac{1}{f(4)}\color{black}\cdot\Big(\frac{1}{f(4)}+\frac{1}{f(5)}+\cdots\Big)&
\end{flalign*}

\begin{flalign*}
&\quad \color{red}-\frac{1}{f(5)}\color{black}\cdot\Big(\frac{1}{f(5)}+\cdots\Big)&
\end{flalign*}

\begin{flalign*}
&\quad \color{red}-\color{black}\cdots&
\end{flalign*}

\noindent which yields the first term in the infinite series,

$$\opM\limits_{i=1}^{} (-1)^{i}\displaystyle\sum{}^{i+1}$$

\noindent or M-series of degree $2$, previously defined as

$$-\Sigma_1 \Sigma_2 = -\sum\limits_{i=1}^\infty \sum\limits_{j=i}^\infty \frac{1}{f(i)} \frac{1}{f(j)}.$$

\noindent The second term of the infinite series

$$\opM\limits_{i=1}^\infty (-1)^{i}\displaystyle\sum{}^{i+1}$$

\noindent is also formed by collecting like terms across the original grouping (again highlighted in \color{red} red \color{black})

\begin{multline*}
$$\displaystyle\sum\limits_{n=1}^\infty \frac{1}{f(n)}\cdot \Big(-\frac{1}{f(1)}\color{red}+\frac{1}{f(1)\cdot f(2)}+\frac{1}{f(1)\cdot f(3)}\color{black}-\frac{1}{f(1)\cdot f(2)\cdot f(3)}\color{red}+\frac{1}{f(1)\cdot f(4)}\color{black}-\frac{1}{f(1)\cdot f(2)\cdot f(4)}-\frac{1}{f(1)\cdot f(3)\cdot f(4)}\\+\frac{1}{f(1)\cdot f(2)\cdot f(3)\cdot f(4)}\color{red}+\frac{1}{f(1)\cdot f(5)}\color{black}-\frac{1}{f(1)\cdot f(2)\cdot f(5)}-\frac{1}{f(1)\cdot f(3)\cdot f(5)}+\frac{1}{f(1)\cdot f(2)\cdot f(3)\cdot f(5)}\\-\frac{1}{f(1)\cdot f(4)\cdot f(5)}+\frac{1}{f(1)\cdot f(2)\cdot f(4)\cdot f(5)}+\frac{1}{f(1)\cdot f(3)\cdot f(4)\cdot f(5)}-\frac{1}{f(1)\cdot f(2)\cdot f(3)\cdot f(4)\cdot f(5)}...\Big)$$
\end{multline*}

\begin{multline*}
$$\quad \displaystyle\sum\limits_{n=2}^\infty \frac{1}{f(n)} \cdot\Big(-\frac{1}{f(2)}\color{red}+\frac{1}{f(2)\cdot f(3)}+\frac{1}{f(2)\cdot f(4)}\color{black}-\frac{1}{f(2)\cdot f(3)\cdot f(4)}\color{red}+\frac{1}{f(2)\cdot f(5)}\color{black}-\frac{1}{f(2)\cdot f(3)\cdot f(5)}\\-\frac{1}{f(2)\cdot f(4)\cdot f(5)}+\frac{1}{f(2)\cdot f(3)\cdot f(4)\cdot f(5)}...\Big)$$
\end{multline*}

\begin{flalign*}
&\quad \displaystyle\sum\limits_{n=3}^\infty \frac{1}{f(n)}\cdot\Big(-\frac{1}{f(3)}\color{red}+\frac{1}{f(3)\cdot f(4)}+\frac{1}{f(3)\cdot f(5)}\color{black}-\frac{1}{f(3)\cdot f(4)\cdot f(5)}...\Big)&
\end{flalign*}

\begin{flalign*}
&\quad \displaystyle\sum\limits_{n=4}^\infty\frac{1}{f(n)}\cdot \Big(-\frac{1}{f(4)}\color{red}+\frac{1}{f(4)\cdot f(5)}\color{black}...\Big)&
\end{flalign*}

\begin{flalign*}
&\quad \displaystyle\sum\limits_{n=5}^\infty \frac{1}{f(n)}\cdot \Big(-\frac{1}{f(5)}...\Big)&
\end{flalign*}

\begin{flalign*}
&\quad \displaystyle\sum\limits_{n=6}^\infty \frac{1}{f(n)}&
\end{flalign*}

\noindent and these terms can be extracted from the original grouping and organized in the same fashion:

\begin{flalign*}
&\quad\Big(\color{red}\frac{1}{f(1)\cdot f(2)}+\frac{1}{f(1)\cdot f(3)}+\frac{1}{f(1)\cdot f(4)}+\frac{1}{f(1)\cdot f(5)}+\cdots\color{black}\Big)\cdot\displaystyle\sum\limits_{n=1}^\infty \frac{1}{f(n)}+&
\end{flalign*}

\begin{flalign*}
&\quad\Big(\color{red}\frac{1}{f(2)\cdot f(3)}+\frac{1}{f(2)\cdot f(4)}+\frac{1}{f(2)\cdot f(5)}+\cdots\color{black}\Big)\cdot\displaystyle\sum\limits_{n=2}^\infty \frac{1}{f(n)}+&
\end{flalign*}

\begin{flalign*}
&\quad\Big(\color{red}\frac{1}{f(3)\cdot f(4)}+\frac{1}{f(3)\cdot f(5)}+\cdots\color{black}\Big)\cdot\displaystyle\sum\limits_{n=3}^\infty \frac{1}{f(n)}+&
\end{flalign*}

\begin{flalign*}
&\quad\Big(\color{red}\frac{1}{f(4)\cdot f(5)}+\cdots\color{black}\Big)\cdot\displaystyle\sum\limits_{n=4}^\infty \frac{1}{f(n)}+\cdots&
\end{flalign*}

\begin{flalign*}
&\quad=&
\end{flalign*}

\begin{flalign*}
&\quad \Big(\color{red}\frac{1}{f(1)\cdot f(2)}+\frac{1}{f(1)\cdot f(3)}+\frac{1}{f(1)\cdot f(4)}+\frac{1}{f(1)\cdot f(5)}+\cdots\color{black}\Big)\cdot\Big(\frac{1}{f(1)}+\frac{1}{f(2)}+\frac{1}{f(3)}+\frac{1}{f(4)}+\frac{1}{f(5)}+\cdots\Big)+&
\end{flalign*}

\begin{flalign*}
&\quad \Big(\color{red}\frac{1}{f(2)\cdot f(3)}+\frac{1}{f(2)\cdot f(4)}+\frac{1}{f(2)\cdot f(5)}+\cdots\color{black}\Big)\cdot\Big(\frac{1}{f(2)}+\frac{1}{f(3)}+\frac{1}{f(4)}+\frac{1}{f(5)}+\cdots\Big)+&
\end{flalign*}

\begin{flalign*}
&\quad \Big(\color{red}\frac{1}{f(3)\cdot f(4)}+\frac{1}{f(3)\cdot f(5)}+\cdots\color{black}\Big)\cdot\Big(\frac{1}{f(3)}+\frac{1}{f(4)}+\frac{1}{f(5)}+\cdots\Big)+&
\end{flalign*}

\begin{flalign*}
&\quad \Big(\color{red}\frac{1}{f(4)\cdot f(5)}+\cdots \color{black}\Big)\cdot\Big(\frac{1}{f(4)}+\frac{1}{f(5)}+\cdots\Big)+\cdots&
\end{flalign*}

\noindent which yields the second term in the infinite series,

$$\opM\limits_{i=2}^{} (-1)^{i}\displaystyle\sum{}^{i+1}$$

\noindent or M-series of degree $3$, previously defined as

$$+\Sigma_1 \Sigma_2 \Sigma_3 = +\sum\limits_{i=1}^\infty \sum\limits_{j=i}^\infty \sum\limits_{k=j+1}^\infty \frac{1}{f(i)} \frac{1}{f(j)} \frac{1}{f(k)}.$$

\noindent The third term of the infinite series

$$\opM\limits_{i=1}^\infty (-1)^{i}\displaystyle\sum{}^{i+1}$$

\noindent is also formed by collecting the next level of terms across the summations in the original grouping

\begin{multline*}
$$\displaystyle\sum\limits_{n=1}^\infty \frac{1}{f(n)}\cdot \Big(-\frac{1}{f(1)}+\frac{1}{f(1)\cdot f(2)}+\frac{1}{f(1)\cdot f(3)}\color{red}-\frac{1}{f(1)\cdot f(2)\cdot f(3)}\color{black}+\frac{1}{f(1)\cdot f(4)}\color{red}-\frac{1}{f(1)\cdot f(2)\cdot f(4)}\\-\color{red}\frac{1}{f(1)\cdot f(3)\cdot f(4)}\color{black}+\frac{1}{f(1)\cdot f(2)\cdot f(3)\cdot f(4)}+\frac{1}{f(1)\cdot f(5)}\color{red}-\frac{1}{f(1)\cdot f(2)\cdot f(5)}-\color{black}\color{red}\frac{1}{f(1)\cdot f(3)\cdot f(5)}\\\color{black}+\frac{1}{f(1)\cdot f(2)\cdot f(3)\cdot f(5)}\color{red}-\frac{1}{f(1)\cdot f(4)\cdot f(5)}\color{black}+\frac{1}{f(1)\cdot f(2)\cdot f(4)\cdot f(5)}+\frac{1}{f(1)\cdot f(3)\cdot f(4)\cdot f(5)}\\-\frac{1}{f(1)\cdot f(2)\cdot f(3)\cdot f(4)\cdot f(5)}...\Big)$$
\end{multline*}

\begin{multline*}
$$\quad \displaystyle\sum\limits_{n=2}^\infty \frac{1}{f(n)} \cdot\Big(-\frac{1}{f(2)}+\frac{1}{f(2)\cdot f(3)}+\frac{1}{f(2)\cdot f(4)}\color{red}-\frac{1}{f(2)\cdot f(3)\cdot f(4)}\color{black}+\frac{1}{f(2)\cdot f(5)}\color{red}-\frac{1}{f(2)\cdot f(3)\cdot f(5)}\\-\color{red}\frac{1}{f(2)\cdot f(4)\cdot f(5)}\color{black}+\frac{1}{f(2)\cdot f(3)\cdot f(4)\cdot f(5)}...\Big)$$
\end{multline*}

\begin{flalign*}
&\quad \displaystyle\sum\limits_{n=3}^\infty \frac{1}{f(n)}\cdot\Big(-\frac{1}{f(3)}+\frac{1}{f(3)\cdot f(4)}+\frac{1}{f(3)\cdot f(5)}\color{red}-\frac{1}{f(3)\cdot f(4)\cdot f(5)}\color{black}...\Big)&
\end{flalign*}

\begin{flalign*}
&\quad \displaystyle\sum\limits_{n=4}^\infty\frac{1}{f(n)}\cdot \Big(-\frac{1}{f(4)}+\frac{1}{f(4)\cdot f(5)}...\Big)&
\end{flalign*}

\begin{flalign*}
&\quad \displaystyle\sum\limits_{n=5}^\infty \frac{1}{f(n)}\cdot \Big(-\frac{1}{f(5)}...\Big)&
\end{flalign*}

\begin{flalign*}
&\quad \displaystyle\sum\limits_{n=6}^\infty \frac{1}{f(n)}&
\end{flalign*}

\noindent which are extracted in the same fashion:

\begin{multline*}
$$\displaystyle\Big(\color{red}-\frac{1}{f(1)\cdot f(2)\cdot f(3)}-\frac{1}{f(1)\cdot f(2)\cdot f(4)}-\frac{1}{f(1)\cdot f(3)\cdot f(4)}-\frac{1}{f(1)\cdot f(2)\cdot f(5)}-\frac{1}{f(1)\cdot f(3)\cdot f(5)}\\-\color{red}\frac{1}{f(1)\cdot f(4)\cdot f(5)}-\cdots\color{black}\Big)\cdot\displaystyle\sum\limits_{n=1}^\infty \frac{1}{f(n)}+$$
\end{multline*}

\begin{flalign*}
&\quad\Big(\color{red}-\frac{1}{f(2)\cdot f(3)\cdot f(4)}-\frac{1}{f(2)\cdot f(3)\cdot f(5)}-\frac{1}{f(2)\cdot f(4)\cdot f(5)}-\cdots\color{black}\Big)\cdot\displaystyle\sum\limits_{n=2}^\infty \frac{1}{f(n)}+&
\end{flalign*}

\begin{flalign*}
&\quad\Big(\color{red}-\frac{1}{f(3)\cdot f(4)\cdot f(5)}-\cdots\color{black}\Big)\cdot\displaystyle\sum\limits_{n=3}^\infty \frac{1}{f(n)}+\cdots&
\end{flalign*}

\begin{flalign*}
&\quad=&
\end{flalign*}

\begin{multline*}
$$\displaystyle\Big(\color{red}-\frac{1}{f(1)\cdot f(2)\cdot f(3)}-\frac{1}{f(1)\cdot f(2)\cdot f(4)}-\frac{1}{f(1)\cdot f(3)\cdot f(4)}-\frac{1}{f(1)\cdot f(2)\cdot f(5)}-\frac{1}{f(1)\cdot f(3)\cdot f(5)}\\\color{red}-\frac{1}{f(1)\cdot f(4)\cdot f(5)}-\cdots\color{black}\Big)\cdot\Big(\frac{1}{f(1)}+\frac{1}{f(2)}+\frac{1}{f(3)}+\frac{1}{f(4)}+\frac{1}{f(5)}+\cdots\Big)+$$
\end{multline*}

\begin{flalign*}
&\quad \Big(\color{red}-\frac{1}{f(2)\cdot f(3)\cdot f(4)}-\frac{1}{f(2)\cdot f(3)\cdot f(5)}-\frac{1}{f(2)\cdot f(4)\cdot f(5)}-\cdots\color{black}\Big)\cdot\Big(\frac{1}{f(2)}+\frac{1}{f(3)}+\frac{1}{f(4)}+\frac{1}{f(5)}+\cdots\Big)+&
\end{flalign*}

\begin{flalign*}
&\quad \Big(\color{red}-\frac{1}{f(3)\cdot f(4)\cdot f(5)}\color{black}\Big)\cdot\Big(\frac{1}{f(3)}+\frac{1}{f(4)}+\frac{1}{f(5)}+\cdots\Big)+\cdots&
\end{flalign*}

\noindent which yields the third term in the infinite series,

$$\opM\limits_{i=3}^{} (-1)^{i}\displaystyle\sum{}^{i+1}$$

\noindent or M-series of degree $4$, previously defined as

$$-\Sigma_1 \Sigma_2 \Sigma_3 \Sigma_4 = -\sum\limits_{i=1}^\infty \sum\limits_{j=i}^\infty \sum\limits_{k=j+1}^\infty \sum\limits_{l=k+1}^\infty \frac{1}{f(i)} \frac{1}{f(j)} \frac{1}{f(k)} \frac{1}{f(l)}.$$

\noindent We now see that a pattern emerges in which it is clear that an infinite sum of infinite sums will be obtained which will add to unity in the numerator of the right-hand side of Eq. 3 if one continues to eliminate terms from the right-hand side of Eq. 4 using the non-sieving method.  It is hypothesized that once all terms are sieved out of the right hand side of the big Zeta function in Eq. 4, a constant will emerge for the M-series in Eq. 3 which takes on a finite value between

$$\color{red}-1\le\color{black} \quad \opM\limits_{i=1}^\infty (-1)^{i}\displaystyle\sum{}^{i+1}\quad \color{red}\le{0}\color{black}$$

\noindent which may reveal something about the frequency of prime numbers yielded by an integer-valued polynomial function $f(n)$ for some limit of the input domain $n \le x$.\\

\noindent Table 3 depicts the M-series values of $M_{x=100}$ and $M_{x=200}$ for the non-sieving application of the Euler zeta function to the domain of the integers $2\le i \le 100$ and $2\le i \le 200$.  It is observed from this table that as the upper limit increases from $x=100$ to $x=200$ that the value of the M-series approaches $-1$.  Table 4 includes the values of $M_{x=100}$ and $M_{x=200}$ for the prime shell function $n^p-(n-1)^p$ for prime powers $2, 3, 5$ and $7$, and it appears that these values will stabilize to some limit as $x$ approaches $\infty$.\\\\\\\\

\begin{center}

\begin{tabular}{|c||c|c|c||c|c|c|}
\hline
prime number $\ge2$ & $\pi \le 100$ & $Ls \le 100$ & ${\displaystyle\opM\limits_{}^{}}_{x=100}$ &  $\pi \le 200$ & $Ls \le 200$ & ${\displaystyle\opM\limits_{}^{}}_{x=200}$ \\
\hline
\textit{i} & 25 & 29.99144 & -0.94812622482360 & 46 & 50.04329 & -0.97060984525939 \\
\hline
\end{tabular}
\end{center}

\begin{center}
\textit{Table 3}\\
\end{center}

\begin{center}

\footnotesize

\begin{tabular}{|c||c|c|c||c|c|c|}

\hline
prime power $p$ & $\Pi \le f(100)$ & $Ms \le f(100)$ & ${\displaystyle\opM\limits_{}^{}}_{x=100}$ &  $\Pi \le f(200)$ & $Ms \le f(200)$ & ${\displaystyle\opM\limits_{}^{}}_{x=200}$ \\
\hline
2 & 44 & 42.75969 & -0.70856869191073 & 76 & 78.48273 & -0.77232394108548 \\
\hline
3 & 43 & 29.01307 & -0.05016737946525 & 72 & 53.06455 & -0.05053523893596 \\
\hline
5 & 18 & 19.71488 & -0.00129463514931 & 32 & 35.92022 &-0.00129463735049 \\
\hline
7 & 24 & 15.71077 & -0.00006682330849 & 40 & 28.56513 & -0.00006682330851 \\
\hline
\end{tabular}
\end{center}

\begin{center}
\textit{Table 4}\\
\end{center}

\noindent It is anticipated that the finite values yielded by the M-series function in non-sieving applications of the Euler zeta function in the case of integer-valued polynomials in general may tell us something about how many prime numbers exist within the ranges that those polynomials were evaluated compared to the number of primes that exist on the real number line less than or equal to some upper limit $x$ that served as the input domain for the base $n$ for those functions.  Thus, the following theorem is proposed:

\begin{theorem}
\noindent When a non-sieving application of the Euler zeta function is applied to process values generated by an integer-valued polynomial, then there is an infinite series

$$\displaystyle\opM\limits_{i=2}^\infty \Big{[(-1)^{i-1}\Big]\displaystyle\sum{}^i}$$

\noindent that arises which adds to unity in the numerator of the product term in the zeta function to make both sides of the equation equal.  The value of this M-series is bounded by $-1$ and $0$ in the non-sieving application of the Euler zeta function.
\end{theorem}

\noindent It is the author's hope that the results of this study will motivate further research into the behavior of prime number frequency among the prime shell function and any other integer-generating polynomial function in general which has the capability to generate prime numbers.\\

\end{document}